\documentstyle[11pt]{article}

\begin{document}

\linespread{1.4}
\newcommand{\Z}{{\Bbb Z}}
\newtheorem{thm}{Theorem}[section]
\newtheorem{defn}{Definition}[section]
\newtheorem{prop}[thm]{Proposition}

\newtheorem{con}[thm]{Conclusion}
\newtheorem{obs}[thm]{Observation}
\newtheorem{cor}[thm]{Corollary}
\newtheorem{lem}[thm]{Lemma}
\newcommand{\pr}{\bf Proof}
\newtheorem{ex}[thm]{Example}

\title { The Classification of Circulant Weighing Matrices of Weight 16 
and Odd Order}
\author {R.\ M.\ Adin\thanks{Supported in part by an internal research grant 
from Bar-Ilan University.}, L.\ Epstein  and  Y.\ Strassler \\
{\it E-mail:}
{\tt radin, epstin and strasler (@macs.biu.ac.il)}\\ 
Department of Mathematics and Computer Science \\
Bar-Ilan University \\ 
Ramat-Gan 52900, Israel}
\date{October 14, 1999}
\maketitle

\begin{abstract}
%\noindent
In this paper we completely classify the circulant weighing matrices of 
weight 16 and odd order. 
It turns out that the order must be an odd multiple of either 21 or 31. 
Up to equivalence, there are two distinct matrices in 
$ CW(31,16) $, one matrix in $CW(21,16)$ and  another one in $CW(63,16)$ 
(not obtainable by Kronecker product from $CW(21,16)$).
The classification uses a multiplier existence theorem. 
	 
\end{abstract}

\normalsize{
\section{Introduction}\label{sec0}
%{\normalsize
In the last
decade we have witnessed an enormous amount of activity in the field of designs. One of the most active subfields is the subject of orthogonal designs. 
It emerged as an attempt to unify the attempts made to close upon the Hadamard 
conjecture which has generated a 
tremendous amount of research in combinatorial matrix theory in the last 
century. 
%It states that for every integer number $n$ divisible by 4 there 
%exists a square matrix $H$ with entries $\pm 1$ and satisfying $HH^{t}=nI$ 
%({\it Hadamard matrix}). with its own subjects and open problems. 
The theory of Hadamard matrices has a lot of applications. For example, in coding theory \cite{sy92}, difference sets \cite{jun92}, spectrometry \cite{hs79}, image processing, image coding, 
	   pattern recognition, sequence filtering \cite{hal79}, genetic algorithms \cite{koe94}, 
	   weighing designs for chemistry and medicine \cite{kal75}. 

One of the tools for investigating Hadamard matrices is circulant weighing
matrices.  
%
%A {\it circulant weighing matrix of order n and weight k} is a square matrix 
%$ W $ of size $ n \times n $ and entries from $ \{-1, \; 0, \; 1 \} $ satisfying~
%$ WW^{t}=kI_{n} $~ with the additional requirement that each row is a right 
%cyclic shift of its predecessor. 
%
Circulant weighing matrices have been known to exist since 1975, when 
A.V. Geramita, J.M. Geramita and J. Seberry \cite{ggs76} observed the existence of a $CW(7,4)$ with first 
row 
	$$ - \ + \ + \ 0 \ + \ 0 \ 0. $$

There are two major classification results of $CW(n,k)$ for fixed weight $k$: 
One by R.~Hain classifying $CW(n,4) $ \cite{hai77}, \cite{eh76}, and the other by 
Y.~Strassler classifying $CW(n,9)$ \cite{str83}, \cite{str}. 

Here is a short list of previously known results about $CW(n,16)$:  

In 1975, Seberry and Whiteman \cite{sw75} proved that 
$CW(q^2+q+1,q^2) \neq \emptyset$ for $q=p^{\alpha}$, $p$ a prime, $\alpha \geq 1$. 
In particular, they constructed one $CW(21,16)$.

In 1980, Eades \cite{ead80} found a $CW(31,16)$ with first row
	$$- \ 0 \ 0 \ 0 \ 0 \ - \ 0 \ + \ 0 \ - \ - \ + \ 0 \ + \ + \ 0 \ 0 \ 0 \ - \ + \ - \ + \ + \ 0 \ 0 \ + \ + \ 0 \ + \ 0 \ 0.$$

In 1995 (published 1998), Strassler \cite{str95} found another CW(31,16), not equivalent to the obtained by Eades.

In this paper we completely classify $CW(n,16)$, for odd values of $n$. In the 
course of study, a new equivalence class was found in $CW(63,16)$.

%In 1996, Muzychuk  proved  existence of a multiplier under certain conditions \cite{muz97}. 
%His results are cornestonea of the current work.
%Note that there are other multiplier theorems of M. Hall and 
%R.L. McFarland but in the current work we will use only Muzychuk's result. 

The paper has the following structure: \\
Definitions and known results appear in the Preliminaries Section (Section \ref{sec1}). 

Section \ref{mul_sec} contains the statement of a multiplier existence 
theorem, including a proof (essentially due to Muzychuk). This result is a 
starting point for the current work. 

The other sections contain various steps of the actual classification process,  
 attempting to find the possible describing sets 
 $P$ and $N$ and the order $n$ of a circulant weighing matrix of weight 16.

      Section \ref{sec2}  introduces orbit length partitions, and 
	   contains a preliminary computation of 
	   all possible pairs of orbit length partitions of $P$ and $N$. 
      
      Section \ref{sec3} shortens the list of possible pairs by using 
	   restrictions  on the number of short orbits. 

      Section \ref{sec4} contains statements and proofs of general lemmas regarding 
      orbit lengths of differences, which  are then used to further reduce the number of possible pairs. 
	       
      Section \ref{sec5} contains a more delicate analysis of the remaining cases, using 
	counting rather than existence arguments.   
      
      Section \ref{sec6} contains final analysis of the few remaining cases, settling conclusively the 
	questions of existence and equivalence.  Computer search is used here. 
      
      The last section (Section \ref{sec7}) contains a summary of the results 
      obtained.

%}

\section{ Preliminaries } \label{sec1}
%{\normalsize
\begin{enumerate}
     \item A {\it Hadamard matrix} $ H=(h_{ij}) $ is a square matrix of order $ n $, with entries $ h_{ij} \in \{-1,1\}$, satisfying $ HH^t=nI $. Occasionally we refer to it as an $ H $ matrix.
     
     \item The {\it Hadamard matrix conjecture} : Hadamard matrices exist for every order $ n $ divisible by 4 \cite{sw75}. 
     
	 The conjecture's status: Still open, although many constructions 
	   of Hadamard matrices are known.
     
	 \item A generalization and construction aid: A {\it Weighing matrix} \- 
	   $ W=~(w_{ij}) $ of {\it order n } and {\it weight k}  is a square 
	   matrix of order $ n $, with entries \- $ w_{ij} \in \{0,1,-1\} $, 
	   satisfying  $$WW^t=~kI_{n}.$$ 
	   $ W(n,k) $ denotes the set of 
	   weighing matrices of order $ n $ and weight $ k $. 
	   We use occasionally ``$ W $ is a $ W(n,k) $'' instead of \- 
	   ``$ W \in W(n,k) $''. \- Also, 
	   occasionally we refer to it as a $ W $ matrix.
     
     \item The {\it weighing matrix conjecture}: Weighing matrices exist for every order $ n $ divisible by 4 and all weights $ 0\leq k \leq n $ \cite{sy92}.
	    The conjecture's status: Still open, although many constructions 
	   of weighing matrices are known.

     \item A basic construction for a weighing matrix $ W(n_{1}n_{2},k_{1}k_{2}) $ is the \-{\it Kronecker product } 
	   of two weighing matrices $ W_{1}(n_{1},k_{1})$, $W_{2}(n_{2},k_{2}) $. 
	   This is the block matrix 
	   $$ W=(w_{ij})=((W_{1})_{ij}W_{2})^{n_{1}}_{i,j=1} .$$
	   We denote this construction by \-$ W=~W_{1} \otimes W_{2} $.   
     
     \item A {\it circulant matrix } is a square matrix in which each row 
	(except the first) is a right cyclic shift of its predecessor. Since 
	the first row of a circulant determines the whole matrix we use 
	the notation $ C=cir(c_{0},c_{1},...,c_{n-1}) $ to denote the 
	circulant matrix 
	  \[ \left( \begin{array}{cccc}
	c_{0}   & c_{1} & \ldots & c_{n-1} \\
	c_{n-1} & c_{0} & \ldots & c_{n-2} \\
	\vdots  \\
	c_{1}   & c_{2} & \ldots & c_{0} \\
	       \end{array}      \right).   \]
	
	 There are many different constructions for $ H $ and $ W $ matrices. Many of them use circulant matrices as construction aids.              
     
     \item A {\it circulant weighing matrix } is a circulant  matrix which is 
	   also a weighing matrix. 
     
     \item Circulant weighing matrices as polynomials:
	Circulant matrices with integer entries form a ring under matrix  addition and 
	multiplication. This ring is 
	isomorphic to the quotient ring \- $R_n:={\bf Z}[x]/~\langle~x^n-1~\rangle$; 
	\- the natural isomorphism
	takes the matrix 
	$$cir(w_0,w_1, \cdots, w_{n-1})$$ 
	into the corresponding {\it Hall polynomial}
	$$w_0+w_1x+ \cdots+ w_{n-1}x^{n-1}.$$ In the ring 
	$R_n$ the weighing property takes the form: 
		$$w(x)w(x^{-1})=~k,$$
	where $x^{-1}:=x^{n-1}$.
	  
     \item  $CW(n,k)$ denotes the set of circulant weighing matrices of order 
	   $n$ and weight $k$. We use occasionally ``$W$ is a $CW(n,k)$'' instead 
	   of \- ``$W~\in~CW(n,k)$''. 
	   
	   Sometimes we identify the matrices in  $CW(n,k)$ with the 
	   corresponding Hall polynomials. 
	   We write \- ``$w(x) \in CW(n,k)$'', where  
	   $w(x)$ is the Hall polynomial of a circulant weighing matrix of 
	   order $n$ and weight $k$.

     \item ${\bf Z}_n$ denotes the ring of integers modulo $n$. 
     
     \item ${\bf Z}^{*}_n$ denotes the multiplicative group of integers modulo $n$, i.e.
	$$ {\bf Z}^{*}_n=\{ t \in {\bf Z}_n \mid gcd(t,n)=1\}.$$
     
     \item \label{eque}
	The circulant weighing matrices $w_1(x), \ w_2(x)$ are {\it 
	equivalent} if they satisfy 
		$$w_2(x)=x^{s}w_1(x^t) \quad  ({\rm in} \ R_n)$$ 
	for some $s \in {\bf Z}_n$, $t \in {\bf Z}^{*}_n$.
     
     \item \label{w-w}
	If $W \in CW(n,k)$ then also \-$-W \in CW(n,k)$. Convention: We refer only to
	one of $W, \ -W$, the one that has more +1's than -1's. 
     
     \item If $CW(n,k) \neq \emptyset$ then $k=s^2$ for some nonnegative integer $s$ \cite{sm76}.
     
     \item \label{mul}
       For $W=cir(w_0,w_1, \cdots ,w_{n-1}) \in CW(n,s^2)$ let \\
       $ P:=\{i \mid w_{i}=1\}$, the {\it positive describing set} of $W$;\\
       $ N:=\{i \mid w_{i}=-1\}$, the {\it negative describing set} of $W$.\\
       Then, using convention (\ref{w-w}) above: 
	$$|P|=\frac{s(s+1)}{2}, \qquad  |N|=\frac{s(s-1)}{2} \quad {\rm \cite{sm76}}.$$
     
     \item The {\it support} of $ w(x) \in CW(n,k)$ is the set
	$$ S=S(w):=\{i \mid w_{i} \neq 0\}.$$
	Hence $ S=P \cup N $. 
     
     \item \label{multip}
	A {\it multiplier} for the  circulant weighing matrix 
	$ w(x) \in CW(n,k) $ is a number $t \in {\bf Z}^{*}_{n}$ 
	such that there exists a shift $ s \in {\bf Z}_{n} $ satisfying  
	$$ w(x)=x^{s}w(x^{t}).$$ 
	If $s=0$ then we say that $t$ is a {\it fixing} multiplier  
	for the circulant $w(x)$. 
	We consider $t=1$ as a {\it trivial} multiplier. From now on, writing 
	``$t$ is a multiplier'' always means ``$t$ is a nontrivial multiplier''.

%     \item \label{muz}
%         The {\it multiplier theorem}: \\
%         Let $w(x) \in CW(n,k)$, $k=s^2$. If $s=p^m$ for a prime number $p$  
%         such that 
%         $gcd(n,p)=1$ (and $m \geq 1$) then $p$ is a multiplier for $w(x)$ \cite{muz97}. 
%       \cite{mcf80}, 
     
     \item \label{str2}
	  Let $w(x) \in CW(n,k)$ have a multiplier $t$ with shift $s$. Then
	$$ tN+s=N \quad {\rm and} \quad tP+s=P. $$
	  In particular, if $t$ fixes $w(x)$ then it also fixes its positive and 
	  negative describing sets: $ tN=N \ {\rm and } \ tP=P$ \cite{str}.

     \item \label{str4}
	 If $w(x) \in CW(n,k)$ with $gcd(n,k)=1$ has a multiplier $t$, then it is 
	 equivalent to some $w'(x) \in CW(n,k)$ for which $t$ is a fixing 
	 multiplier \cite{str}.

     \item Let $ t~\in~{\bf Z}^{*}_{n}. \ {\rm A \ subset} \ Z \subset {\bf Z}_{n} $ 
	is called a {\it t-orbit} if there exists an element 
	$ z \in {\bf Z}_{n} $ such that 
	$$ Z=\{ t^{i}z~\pmod{n}~\mid i \in \bf Z\} .$$
	Denote by $ol(z)$ the {\it orbit length} of $z\in{\bf Z}_{n}$, i.e., the number 
	of elements in the $t$-orbit containing $z$.

    \item A {\it multiset} is a ``set'' in which repetitions of elements are 
	allowed. We distinguish it from a regular set by using brackets [] instead of braces \{\}. 
	For example, 
	$$X=[a,a,b,c,c,c]=[a^2,b,c^3].$$
    
    \item Two multisets can be merged together to form a new multiset by the 
	{\it adjunction} ($\&$) operation, which is the union with repetitions 
	counted. 
	For example,
	   $$X=[a^2,b,c^3], \ Y=[a,b^5,d]$$
	   $$X \& Y=[a^3,b^6,c^3,d].$$

    \item For a multiset $X \subseteq {\bf Z}_{n}$, 
	  $$\triangle X:=[x_1-x_2 \mid \ x_1,x_2 \in X, \  x_1 \neq x_2].$$
    
    \item For multisets $X, \ Y \subseteq {\bf Z}_{n}$, 
	  $$\overline{\triangle} X,Y:= [\pm(x-y) \mid \ x \in X, \ y \in Y].$$
    
    \item \label{eq_p_n}
	The CW multiset equation: If 
	$P$ and $N$ are the (positive and negative) describing sets of some
	$W \in CW(n,k)$ then 
	$$\triangle P \& \triangle N= \overline{\triangle}P,N \quad \cite{str}.$$
	Note that $P$ and $N$ are sets, but 
	$\triangle P$, $\triangle N$ and $ \overline{\triangle}P,N$ 
	are (in general) multisets.

   \item  \label{adjpas} 
	If $q=2^t$ and $i$ is even then \-
	$CW(\frac{q^{i+1}-1}{q-1},q^i) \neq \emptyset$ \- \cite{adjp95}, 
	 \cite{as96}.

\end{enumerate}

%}

\section{A Multiplier Existence Theorem} \label{mul_sec}

A fundamental result, on which the current classification is based, is the following 
multiplier existence theorem. This is by now a folklore result, quoted 
(and sometimes reproved) in many sources - e.g., Mcfarland \cite{mcf80}, 
Lander \cite{lan83}, Arasu \cite{as98}, Jungnickel \cite{jun92}, 
Muzychuk \cite{muz97}.
In order to make the paper self-contained, we include here a relatively short 
and elegant proof, basically due to Muzychuk \cite{muz97}, and adapted to 
suit our context. 

\begin{thm}  
	 Let $w(x) \in CW(n,k)$, $k=s^2$. If $s=p^m$ for a prime $p$  
	 such that 
	 $gcd(n,p)=1$ (and $m \geq 1$), then $p$ is a multiplier for $w(x)$.
\end{thm}
\begin{pr}
{\normalsize { \ }\\
  Let $b$ be the maximal nonnegative integer such that       
	$$ w(x^{p^j})w(x^{-1}) \equiv 0 \pmod{p^b} \quad ({\rm in} \ R_n),$$ 
  for all nonnegative integers $j$. Of course, since 
	$$ w(x)w(x^{-1})=k=p^{2m}  \quad ({\rm in} \ R_n),$$ 
  necessarily $b \leq 2m$. We shall show that $b=2m$. 

  Indeed, let $j_0$ be a nonnegative integer such that 
	$$ w(x^{p^{j_0}})w(x^{-1})=p^b v_0(x)$$ 
  for some $v_0(x) \in R_n$ such that $v_0(x) \not\equiv 0 \pmod{p}$. 
  Then 
   $$ p^{2b} v_0(x^{p^{j_0}}) v_0(x)=
       w(x^{p^{2j_0}}) w(x^{-p^{j_0}}) w(x^{p^{j_0}}) w(x^{-1})=
       w(x^{p^{2j_0}}) w(x^{-1}) p^{2m} \equiv 0 \pmod{p^{b+2m}},$$
  by the definition of $b$. If $b<2m$ then it follows that 
	$$v_0(x^{p^{j_0}}) v_0(x) \equiv 0 \pmod{p^{2m-b}} $$ 
  so that, in particular (since $v_0(x^p) \equiv v_0(x)^p \pmod {p}$): 
	$$v_0(x)^{p^{{j_0}+1}} \equiv v_0(x^{p^{j_0}}) v_0(x) \pmod{p} \equiv 0 \pmod{p}. $$          
  It follows that $v_0(x)$ is a nilpotent element in the ring  
  ${\bf Z}_p [x] / \langle x^n-1 \rangle$, which is the group algebra over  
  ${\bf Z}_p$ of the cyclic group of order $n$, and is therefore semisimple 
  (since $p \not| \; n$). Thus $v_0(x) \equiv 0 \pmod{p}$, a contradiction. 
  We have shown that $b=2m$, and in particular (for $j=1$): 
	 $$ w(x^p)w(x^{-1}) \equiv 0 \pmod {p^{2m}}  \quad ({\rm in} \ R_n).$$  
  Let $v(x) \in R_n$ satisfy 
	 $$ w(x^p)w(x^{-1})=p^{2m}v(x)  \quad ({\rm in} \ R_n).$$ 
  Then 
   $$ p^{4m} v(x) v(x^{-1})=
       w(x^{p}) w(x^{-1}) w(x^{-p}) w(x)= 
       w(x^{p}) w(x^{-p}) w(x) w(x^{-1})=p^{2m}p^{2m}$$  
  so that 
	$$ v(x)v(x^{-1})=1  \quad ({\rm in} \ R_n).$$ 
  Computing the coefficient of $x^0=1$ on both sides of the equation, we 
  conclude that if 
	$$v(x)=\sum_{i=0}^{n-1}v_i x^i  \quad (v_i \in {\bf Z})$$
  then 
	$$\sum_{i=0}^{n-1}v_i^2 =1 \quad ({\rm in} \ {\bf Z}). $$ 
  Thus exactly one of the $v_i$ is nonzero (and equal to $\pm 1$), so that 
       $$ w(x^p)w(x^{-1})=p^{2m}v(x)=\pm p^{2m} x^i  \quad ({\rm in} \ R_n)$$  
  for some $0 \leq i \leq n-1$. Thus
     $$\pm p^{2m} x^i w(x)= w(x^p)w(x^{-1})w(x)= w(x^p) p^{2m} \quad ({\rm in} \ R_n),$$ 
     $$\pm x^i w(x)= w(x^p) \quad ({\rm in} \ R_n).$$ 
  Obviously, the ``$\pm$" is actually ``$+$" (e.g., since the sum of the 
  coefficients of $w(x)$ is nonzero by (\ref{mul}) from Section 2). Therefore $p$ 
  is a multiplier for $w(x)$. 
			\begin{flushright}
				$\diamondsuit$
			\end{flushright}
}
\end{pr}

\section{ Orbit-Length Partitions } \label{sec2}

%{\normalsize
The present work concerns $ CW(n,16)$ where $n$ is odd.
Here $k=s^{2}=16, \ s=4$, and $gcd(n,2)=1$. 
Thus, according to the multiplier existence theorem (Theorem 3.1), 
$t=2$ is a multiplier for each $w(x) \in CW(n,16)$.  
By (\ref{mul}) in Section 2, 
$$ |P|=\frac{s(s+1)}{2}=\frac{4(4+1)}{2}=10, \qquad  
|N|=\frac{s(s-1)}{2}=\frac{4(4-1)}{2}=6.$$ 
By claim (\ref{str4}) in Section 2 we can assume, without  
loss of generality, that $t=2$ is a fixing multiplier for $w(x)$.  
The sets $P, \ N \subseteq {\bf Z}_{n}$ are then closed
under multiplication by the multiplier $t \in  {\bf Z}^{*}_{n} $ (by 
claim (\ref{str2})). It follows that $P$ and $N$ are unions of $t$-orbits.  

 Write now the $t$-orbits 
within $P$ in order of increasing length, that is: 
$ P = {C_{1}} \cup {\cdots} \cup {C_{m}} $, where $ m $ is the number of $t$-orbits 
in $P$, and $ |C_{i}| \leq |C_{i+1}| \; (\forall i)$. 
Denote $ l_{i}=|C_{i}| $ and obtain the 
{\it orbit length partition} of $ P$:  
	      $$ olp(P)=(l_{1},\cdots,l_{m}) \quad (l_{1} \leq \cdots \leq l_{m}).$$ 
If $l_{i}=l_{i+1}= \cdots =l_{i+d_{i}-1}$, write $l^{d_{i}}_{i}$ instead of
$d_{i}$ times $l_{i}$. We shall sometimes use the shorter notation 
$$olp(P)=l^{d_1}_{1}l^{d_2}_{2} \cdots l^{d_k}_{k} \quad (k \leq m)$$
instead of 
$$olp(P)=(l^{d_1}_{1},l^{d_2}_{2}, \cdots, l^{d_k}_{k}) \quad (k \leq m).$$

Define $ olp(N)$ in a similar way.

%$ N = {B_{1}} \cup {\cdots} \cup {B_{h}} $, where $h$ is the number of $t$-orbits 
%in $N$, and \- $|B_{i}|~\leq~|B_{i+1}| \; (\forall i)$. 
%Denote $ f_{i}=|B_{i}| $ and thus obtain 
%the {\it orbit length partition} of $ N$:  
%         $$olp(N)=(f_{1},\cdots,f_{h}) \quad (f_{1} \leq \cdots \leq f_{h}).$$
%Again, if $f_{i}=f_{i+1}= \cdots =f_{i+d_{i}-1}$ write them as $f^{d_{i}}_{i}$ 
%and sometimes write 
%$$olp(N)=f^{d_1}_{1}f^{d_2}_{2} \cdots f^{d_g}_{g} \quad (g \leq h)$$
%instead of 
%$$olp(N)=(f^{d_1}_{1},f^{d_2}_{2}, \cdots, f^{d_g}_{g}) \quad (g \leq h).$$

\begin{ex}

$$w(x)=-x+x^2-x^3+x^5+x^6+x^7+x^8-x^9+x^{11} \in CW(13,9).$$ 
%{\normalsize
{\rm This weighing circulant has $ t=3$ as a multiplier.} \\
$ |N|=3,\ |P|=6.\\
N=\{1,3,9\}, \   olp(N)=3^1;\\
 P=\{2,6,5,8,11,7\}=C_1 \cup C_2, \ {\rm where} \ C_1=\{2,6,5\} \ {\rm and} \ 
 C_2=\{8,11,7\}. \rm \ Hence \   olp(P)=3^2$.
%
%       \begin{flushright}
%               $\diamondsuit$
%       \end{flushright}
\end{ex}

Let us start by listing all possible partitions of $ N$ with  
$|N|=6 $ and  of $P$ with $|P|=10$:\\

$ olp(N) \in \{ 1^{6}, 1^{4}2^{1}, 1^{2}2^{2}, 2^{3}, 1^{3}3^{1}, 1^{1}2^{1}3^{1}, 3^{2},  
   1^{2}4^{1}, 2^{1}4^{1}, 1^{1}5^{1}, 6^{1} \}. $\\

$ olp(P) \in  \{ 1^{10}, 1^{8}2^{1}, 1^{6}2^{2},1^{4}2^{3}, 1^{2}2^{4}, 2^{5}, 1^{7}3^{1}, 1^{5}2^{1}3^{1}, 1^{3}2^{2}3^{1}, 1^{1}2^{3}3^{1}, 1^{4}3^{2}, 
1^{2}2^{1}3^{2}, 2^{2}3^{2}, 1^{1}3^{3}, \\ 1^{6}4^{1}, 1^{4}2^{1}4^{1},1^{2}2^{2}4^{1}, 2^{3}4^{1},1^{3}3^{1}4^{1}, 1^{1}2^{1}3^{1}4^{1}, 3^{2}4^{1},1^{2}4^{2},2^{1}4^{2},
 1^{5}5^{1}, 1^{3}2^{1}5^{1}, 1^{1}2^{2}5^{1}, 1^{2}3^{1}5^{1}, 2^{1}3^{1}5^{1}, \\ 1^{1}4^{1}5^{1}, 5^{2}, 1^{4}6^{1}, 1^{2}2^{1}6^{1}, 2^{2}6^{1}, 
1^{1}3^{1}6^{1}, 4^{1}6^{1},
1^{3}7^{1}, 1^{1}2^{1}7^{1}, 3^{1}7^{1},1^{2}8^{1}, 2^{1}8^{1}, 1^{1}9^{1}, 10^{1} \}.$

\section{ The Number of Short Orbits } \label{sec3}

Each of the $t$-orbits ($t=2$) encountered in the current classification has 
length at most 10. Therefore, if we determine for each $1 \leq i \leq 10$ 
the number of orbits of length $i$ in ${\bf Z}_{n}$, we shall be able to 
exclude orbit length partitions requiring more than this number of orbits, 
and thus reduce the size of our search space. It turns out that $1 \leq i \leq 3$ suffice 
for the basic elimination process, but the cases $4 \leq i \leq 6$ will also be 
needed in Section 7. 

In general, an element $a \in {\bf Z}_{n}$ has orbit length dividing $i$ iff
	    $$ 2^ia \equiv a \pmod {n} ;$$
	    $$ (2^{i}-1)a \equiv 0 \pmod {n} ;$$                                
	    $$ a=\frac{nk}{2^{i}-1} \quad (k \in \{0,\cdots,2^{i}-2\}). $$
Of course, to conclude that the orbit length is exactly $i$, one has to exclude all 
the proper divisors of $i$. This is easy when $i$ is 1 or a prime number, and is 
not too difficult for other small values of $i$. We are interested in the cases 
$1 \leq i \leq 6$. 

\begin{enumerate}
      \item $ \underline{i=1} $:\\ 
	    Let $ a \in {\bf Z}_{n} $ have $ ol(a)=1 $. Then:\\
	    $$ 2^1a \equiv a \pmod {n} ;$$
	    $$ a \equiv 0 \pmod {n} .$$                                
	    Thus there is exactly one element in $ {\bf Z}_{n} $, namely 0, 
	    whose orbit length is equal to 1.
      \item $ \underline{i=2} $:\\ 
	    Let $ a \in {\bf Z}_{n} $ have $ ol(a)=2 $. Then:\\
	    $$ 2^2a \equiv a\pmod {n} ;$$
	    $$ 3a \equiv 0 \pmod {n} ;$$                                
	    $$ a=\frac{nk}{3} \quad (k \in \{0,1,2\}). $$
	    The value $ k=0 $ is impossible, since then $ a=0 $ and   
	    $ol(a)=1 $. Thus there are at most two elements in $ {\bf Z}_{n} $ with 
	    orbit length equal to 2. They form a single orbit:\\
	    $$ (\frac{n}{3},\frac{2n}{3}) .$$
	    This orbit exists iff $ n $ is divisible by 3.
      \item $ \underline{i=3} $:\\ 
	    Let $ a \in {\bf Z}_{n} $ have $ ol(a)=3 $. Then:\\
	    $$ 2^3a \equiv a \pmod {n} ;$$
	    $$ 7a \equiv 0 \pmod {n} ;$$                                
	    $$ a=\frac{nk}{7} \quad (k\in \{0,\cdots,6\}). $$
	    Again $ k=0 $ is impossible.
	    Thus there are at most six elements in $ {\bf Z}_{n} $ with orbit length            
	    equal to 3. It follows that there are at most two orbits of length 3:\\
	    $$ (\frac{n}{7},\frac{2n}{7},\frac{4n}{7}) ;$$
	    $$ (\frac{3n}{7},\frac{6n}{7},\frac{5n}{7}) .$$
	    Each of these orbits exists iff $ n $ is divisible by 7.
    \item $ \underline{i=4} $:\\ 
	The orbit length of $a \in {\bf Z}_n$ divides 4 iff  
	    $$ a=\frac{nk}{15} \quad (k \in \{0,\cdots,14\}). $$
	  The cases $k \in \{0,5,10\}$ lead to shorter orbits (length 1 or 2). The cases 
	  $k \in \{3,6,9,12\}$ are possible whenever $5 \mid n $. Thus there is at least one 
	  orbit of length 4 
	  in ${\bf Z}_n$ iff $n$ is divisible by 5, and there are three different 
	  orbits of length 4 iff $n$ is divisible by 15. 
 
      \item $\underline{ i=5} $:\\ 
	    The orbit length of $ a \in {\bf Z}_{n} $ divides 5 iff\\
	    $$ 2^5a \equiv a \pmod {n} ;$$
	    $$ 31a \equiv 0 \pmod {n} ;$$                                
	    $$ a=\frac{nk}{31} \quad (k \in \{0,\cdots,30\}). $$
	    The value $ k=0 $ leads to $a=0$ with orbit length 1. 
	    Thus there are at most 30 elements in $ {\bf Z}_{n} $ 
	    with orbit length equal to 5. In other words, there are at most 
	    $ 30 \div 5=6 $ orbits of length 5. 
	    Each of them exists iff $n$ is divisible by 31.
      \item $ \underline{i=6} $:\\ 
	    The orbit length of $ a \in {\bf Z}_{n} $ divides 6 iff
	    $$ 2^6a \equiv a \pmod {n} ;$$
	    %$$ 63a \equiv 0 \pmod {n} ;$$                                
	    $$ a=\frac{nk}{63} \quad (k \in \{0,\cdots,62\}). $$ 
	    We have to exclude orbit length 1,2 and 3 (the proper divisors of 6).  
	    Note that $63=3 \times 3 \times 7$. Using the analysis of previous 
	    cases, we get: 

    \begin{itemize}
	\item $ ol(a)=1 $ iff $ k=0 $.
	\item $ ol(a)=2 $ iff $ 21 | k $ and $k \neq 0$, i.e. $k \in \{21,42\}$. 
	\item $ ol(a)=3 $ iff $ 9 | k $ and $k \neq 0$, i.e. $k \in \{9,18,27,36,45,54\}$. 
     \end{itemize}

	We are left with at most $ 63-(1+2+6)=54 $ elements in $ {\bf Z}_{n} $ with orbit length equal to 6. 
	It follows that there are at most $54 \div 6=9$ orbits of length 6.
%
%    $$ (\frac{n}{63},\frac{2n}{63},\frac{4n}{63},\frac{8n}{63},\frac{16n}{63},\frac{32n}{63}), $$
%    $$ (\frac{3n}{63},\frac{6n}{63},\frac{12n}{63},\frac{24n}{63},\frac{48n}{63},\frac{33n}{63}), $$
%    $$ (\frac{5n}{63},\frac{10n}{63},\frac{20n}{63},\frac{40n}{63},\frac{17n}{63},\frac{34n}{63}), $$
%    $$ (\frac{7n}{63},\frac{14n}{63},\frac{28n}{63},\frac{56n}{63},\frac{49n}{63},\frac{35n}{63}), $$
%    $$ (\frac{11n}{63},\frac{22n}{63},\frac{44n}{63},\frac{25n}{63},\frac{50n}{63},\frac{37n}{63}), $$
%    $$ (\frac{13n}{63},\frac{26n}{63},\frac{52n}{63},\frac{41n}{63},\frac{19n}{63},\frac{38n}{63}), $$
%    $$ (\frac{15n}{63},\frac{30n}{63},\frac{60n}{63},\frac{57n}{63},\frac{51n}{63},\frac{39n}{63}), $$
%    $$ (\frac{23n}{63},\frac{46n}{63},\frac{29n}{63},\frac{58n}{63},\frac{53n}{63},\frac{43n}{63}), $$
%    $$ (\frac{31n}{63},\frac{62n}{63},\frac{61n}{63},\frac{59n}{63},\frac{55n}{63},\frac{47n}{63}), $$

    Some of these orbits exist even if $63 \not | \; n$. Indeed, 
    if $n$ is divisible by 63 then there are 9 orbits; if $n$ is  
    divisible by 21 but not by 63 then there are two orbits ($3 | k $ but 
    $9 \not| \; k$ and $21 \not| \; k$, i.e., $k \in \{3,6,12,15,24,30,33,39,48,51,57,60\}$); and if $n$ is divisible
    by 9 but not by 63 then there is only one orbit ($7 | k $ but $21 \not| \; k$, i.e., $k \in \{7,14,28,35,49,56\}$) .
 
\end{enumerate}

%\normalsize{
    
%    In a similar way we can check that the expecting (by previos section)     
%    number of different orbits of length $i$ exist for  $5 \leq i \leq 10$  
%    under certain conditions on $n$.
%    Thus we must choose the order of the circulant weighing matrix accordingly
%to the given $olp(P), \ olp(N)$ (in other words, accordingly to the number and the length of different orbits in $P$, $N$).\\

Let us return now to $olp(P)$ and $olp(N)$. From the above analysis 
of the cases $i \in \{1,2,3\}$  it follows that we can delete
from our list all orbit length partitions that contain 
$$ 1^{a} \ {\rm with}\ a>1, \ 2^{a}  \  {\rm with} \ a>1, \ {\rm or} \ 3^{a} \ {\rm with} \ a>2. $$

Thus there remain only the following orbit length partitions :

%}

\begin{flushleft}
$ olp(N) \in \small{\{ 1^{1}2^{1}3^{1},\;  3^{2},\;  2^{1}4^{1},\; 1^{1}5^{1},\; 6^{1} \}}. $

$ olp(P) \in \small{\{ 1^{1}2^{1}3^{1}4^{1},\; 3^{2}4^{1},\; 2^{1}4^{2},\; 2^{1}3^{1}5^{1},\;
   1^{1}4^{1}5^{1},\; 5^{2},\; 1^{1}3^{1}6^{1},\; 4^{1}6^{1}, \;1^{1}2^{1}7^{1},\; }$ \\
   $\small{3^{1}7^{1},\; 2^{1}8^{1},\; 1^{1}9^{1},\;  10^{1} \}}.$

\end{flushleft}

%\normalsize{ 
Overall, we still have $ 5 \times 13=65 $ cases to check.
Actually, the restrictions on the number of orbits of given size apply not only
to each of $olp(P)$ and $olp(N)$ separately, but to the combined partition $olp(P) \cup olp(N)$ as well.

Applying this condition, there remain only 41 possible pairs $(olp(P), olp(N))$: 

\large{\bf Table 1}: \normalsize{
Initial list of pairs of orbit length partitions}

\begin{center}

\begin{tabular}{||l|l|l||}      \hline
  $ \# $ & $olp(P)$               & $olp(N)$           \\ \cline{1-3}

     1   & $5^{2}$                & $1^{1}2^{1}3^{1}$   \\ \hline 
     2   & $4^{1}6^{1}$           & $1^{1}2^{1}3^{1}$   \\ \hline 
     3   & $3^{1}7^{1}$           & $1^{1}2^{1}3^{1}$   \\ \hline 
     4   & $10^{1}$               & $1^{1}2^{1}3^{1}$   \\ \hline 
     5   & $2^{1}4^{2}$           & $3^{2}$             \\ \hline 
     6   & $1^{1}4^{1}5^{1}$      & $3^{2}$             \\ \hline

\end{tabular}

\end{center}
     
\begin{center}

\begin{tabular}{||l|l|l||}      \hline
  $ \# $ & $olp(P)$               & $olp(N)$           \\ \cline{1-3}

     7   & $5^{2}$                & $3^{2}$             \\ \hline 
     8   & $4^{1}6^{1}$           & $3^{2}$             \\ \hline 
     9   & $1^{1}2^{1}7^{1}$      & $3^{2}$             \\ \hline 
     10  & $2^{1}8^{1}$           & $3^{2}$             \\ \hline 
     11  & $1^{1}9^{1}$           & $3^{2}$              \\ \hline 
     12  & $10^{1}$              & $3^{2}$             \\ \hline 
     13  & $3^{2}4^{1}$          & $2^{1}4^{1}$        \\ \hline 
     14  & $1^{1}4^{1}5^{1}$     & $2^{1}4^{1}$        \\ \hline 
     15  & $5^{2}$               & $2^{1}4^{1}$        \\ \hline 
     16  & $1^{1}3^{1}6^{1}$     & $2^{1}4^{1}$        \\ \hline 
     17  & $4^{1}6^{1}$          & $2^{1}4^{1}$        \\ \hline 
     18  & $3^{1}7^{1}$          & $2^{1}4^{1}$        \\ \hline 
     19  & $1^{1}9^{1}$          & $2^{1}4^{1}$        \\ \hline 
     20  & $10^{1}$              & $2^{1}4^{1}$        \\ \hline 
     21  &$ 3^{2}4^{1} $         & $1^{1}5^{1}$        \\ \hline 
     22  & $2^{1}4^{2} $         & $1^{1}5^{1}$        \\ \hline 
     23  & $2^{1}3^{1}5^{1}$     & $1^{1}5^{1}$        \\ \hline  
     24  &$ 5^{2}$               & $1^{1}5^{1}$        \\ \hline 
     25  & $4^{1}6^{1} $         & $1^{1}5^{1}$         \\ \hline 
     26  & $3^{1}7^{1}$          & $1^{1}5^{1}$         \\ \hline 
     27  & $2^{1}8^{1}$          & $1^{1}5^{1}$         \\ \hline 
     28  & $10^{1}$              & $1^{1}5^{1}$        \\ \hline 
     29  & $1^{1}2^{1}3^{1}4^{1}$& $6^{1}$             \\ \hline 
     30  &$ 3^{2}4^{1}$          & $6^{1}$         \\ \hline 
     31  & $2^{1}4^{2}$          & $6^{1}$         \\ \hline 
     32  & $2^{1}3^{1}5^{1}$     & $6^{1}$         \\ \hline 
     33  & $1^{1}4^{1}5^{1}$     & $6^{1}$         \\ \hline 
     34  & $5^{2}$               & $6^{1}$         \\ \hline 
     35  & $1^{1}3^{1}6^{1}$     & $6^{1}$         \\ \hline   
     36  & $4^{1}6^{1}$          & $6^{1}$         \\ \hline
     37  & $1^{1}2^{1}7^{1}$     & $6^{1}$         \\ \hline
     38  & $3^{1}7^{1}$          & $6^{1}$         \\ \hline

\end{tabular}

\end{center}
     
\begin{center}

\begin{tabular}{||l|l|l||}      \hline
  $ \# $ & $olp(P)$              & $olp(N)$           \\ \cline{1-3}

     39  & $2^{1}8^{1}$          & $6^{1}$         \\ \hline
     40  & $1^{1}9^{1}$          & $6^{1}$         \\ \hline            
     41  & $10^{1}$              & $6^{1}$         \\ \hline

\end{tabular}

\end{center}
{ \ }\\

\section{ Auxiliary Lemmas on Differences }\label{sec4}

In this section we shall formulate and prove lemmas, concerning orbit lengths of differences, 
that will be useful in eliminating more cases. 
Recall the notation $ol(a)$ for the orbit length of $a \in {\bf Z}_{n}$. 
%}
%{\ }\\
\begin{obs}\label{cl0} 
%\normalsize{
	If $ol(a)=1$ and $ol(b)=k>1$, then $ol(a-b)=k$.
%               }
\end{obs}
%\normalsize{
	Indeed, by the previous section $a=0$, hence 
	   $ol(a-b)=ol(-b)=ol(b)=k.$
%          }
%
%{\ }\\

\begin{obs}\label{cl01} %\normalsize{
	If $a \neq b$, then $ol(a-b) > 1$.
		%}
\end{obs}

%{ \ }\\

Denote by $ gcd(a,b)$ the {\it greatest common divisor} and by 
$lcm(a,b)$ the {\it least common multiple} of the integers $a$ and $b$. 

\begin{lem}\label{cl1} %\normalsize{
	If $ ol(a)=k $, $ ol(b)=l $ and $ ol(a-b)=m $ then
	\begin{enumerate}
	      \item $ m \mid lcm(k,l); $  
	      \item $ k \mid lcm(l,m); $   
	      \item $ l \mid lcm(m,k); $   
	\end{enumerate}     %}
\end{lem}                                                   
\begin{pr}: 
{\rm
	  \begin{enumerate}
	      \item Since $t^{k}a=a$ and $t^{l}b=b$, it follows that
		    $t^{lcm(k,l)}a=a$ and $t^{lcm(k,l)}b=b$ so that 
		    $t^{lcm(k,l)}(a-b)=a-b$ as well.
		    Since, for $i \geq 0$,
		    $t^{i}(a-b)=a-b$ iff $ol(a-b) \mid i$, it follows
		    that $m=ol(a-b)$ divides $lcm(k,l)$.
	      
	      \item Since $a=b-(b-a)$ and $ol(b-a)=ol(a-b)=m$, replacing 
		    $a, \ b, \ a-b $ by $b, \ b-a, \ a$, respectively, gives 
		    $ k \mid lcm(l,m) $.
	      \item Similarly, since $b=a-(a-b)$, replacing  
		    $a, \ b, \ a-b $ from the first case by $a, \ a-b, \ b$, 
		    respectively, gives $ l \mid lcm(m,k) $.
	\end{enumerate}
		       }
	\begin{flushright}
		$\diamondsuit$
	\end{flushright}
\end{pr}

\begin{cor}\label{cl11} %\normalsize{
	If $ ol(a)=ol(b)=k, \ a \neq b $, and $k$ is prime, then $ol(a-b)=k$.
		     %}
\end{cor}

{\rm
Indeed, by Lemma \ref{cl1}, $ol(a-b) \mid k$. Since $k$ is prime, 
$ol(a-b)=1$ or $ol(a-b)=k$. The first option is impossible because $ a \neq b$;
hence $ol(a-b)=k$.
	  }
%{\ }\\

\begin{lem} \label{cl2} %\normalsize{
	If $ ol(a)=k, \ ol(b)=l $ and $ gcd(k,l)=1 $, then $ ol(a-b)=kl$. %}
\end{lem}
\begin{pr}:{\rm  \\
		Let $ m:=ol(a-b) $.
	\begin{itemize}
		   \item By Lemma \ref{cl1}, $ m \mid lcm(k,l) $; but $gcd(k,l)=1$,
		    so that $ lcm(k,l)=kl $. Hence $ m=k'l'$, where $ k' \mid k $ and 
		    $ l' \mid l $. We shall prove that  $ k'=k$ and $ l'=l $.
	      \item $ l \mid lcm(k,m) $, by Lemma \ref{cl1}. Notice that 
		    $ gcd(k,l)=1 $, hence $ l \mid m $. Thus $l'=l$.
	      \item $ k \mid lcm(m,l) $, by Lemma \ref{cl1}. In other words,
		    $ k \mid lcm(k'l,l)=k'l $. Notice that $ gcd(k,l)=1 $, 
		    hence $ k \mid k' $. Therefor $ k'=k $ and $ m=kl $.
	\end{itemize}
		 }
	\begin{flushright}
		$\diamondsuit$
	\end{flushright}
\end{pr}
%{\ }\\

\begin{lem}\label{cl4} %\normalsize{
       Let $ ol(a)=k $, $ ol(b)=m $, and $ ol(a-b)=l $. If $m$ is prime      
       then exactly one of the following holds:
       \begin{enumerate}
	    \item $l=km, \quad m \not| \; k$.  
	    \item $ l=k, \quad m \mid k $. 
	    \item $ l=\frac{k}{m}, \quad m \mid k, \quad m \not| \; l $.
       \end{enumerate}                                        
		 % }
\end{lem}
\begin{pr}: {\rm
It is clear that the three cases are mutually exclusive. Thus it suffices to show that at 
least one of them holds. 
 
 If $ m \not| \; k$ then $gcd(k,m)=1$, hence by Lemma \ref{cl2} $l=km$. This gives case 1. 
      
 If $m \mid k$  then one of the following holds: 
 Either $ m \mid l$, and then $ lcm(k,m)=k$ and $ lcm(l,m)=l$. Thus, by Lemma~\ref{cl1}, 
 $ l \mid k$ and $ k \mid l$ and therefore $l=k$. This is case 2. 

 Alternatively, $ m \not| \; l $. Then $ lcm(k,m)=k$ and $ lcm(l,m)=lm$, so that, by Lemma~
	\ref{cl1}, $ l \mid k$ and $ k \mid lm$.
	Writing $ k=k'm$, we get
	$ l \mid k'm$ and $ k'm \mid lm$. Hence
	$ l \mid k'$ and $ k' \mid l$. Thus $ l=k'=\frac{k}{m}$,
	$ m \mid k$ but $m \not| \; l$. This is case 3. 
	    
		     }
	\begin{flushright}
		$\diamondsuit$
	\end{flushright}
\end{pr}

%{\ }\\

\begin{lem}\label{cl6} %\normalsize{
Let $ ol(a)=k, \ ol(b)=m $, and $ol(a-b)=l $. If $ k=k'u $ and $ m=m'u$ with  
$gcd(k',m')=1$ and $u$ prime
then:
	\begin{enumerate}
	  \item Either $ l=k'm' \  or \  l=uk'm' $.
	  \item If either $ u \mid k' \  or \ u \mid m' $, then $ l=uk'm' $.  
	\end{enumerate}
		 %}
\end{lem}
\begin{pr}: \\ {\rm 
By Lemma \ref{cl1}, 
     $$ k \mid lcm(l,m) \; {\Rightarrow} \; k'u \mid lcm(l,m'u) \; \stackrel{gcd(k', m')=1}{\Longrightarrow} \; k' \mid l $$
     and
     $$ {m \mid lcm(l,k)} \; \Rightarrow \; {m'u \mid lcm(l,k'u) } \; \stackrel{gcd(k', m')=1}{\Longrightarrow} \; {m' \mid l} .$$
     Hence $ k'm' \mid l $.
\begin {enumerate}
   \item 
	 Since $ lcm(k,m)=uk'm' $, it follows from Lemma \ref{cl1} that 
	 $ l~\mid~uk'm' $. \- Since also $k'm' \mid l$ and $u$ is prime, it follows that 
	 either $ l=k'm'\  {\rm or} \ l=uk'm'$.
   \item Assume, for example, that $ u \mid k'$.\\  
     Write $ k'=uk''$. By Lemma \ref{cl1}, $ {k \mid lcm(l,m)} $. If $l \neq uk'm'$  
     then 
		$$ l=k'm'=uk''m'=k''m $$ 
     and we obtain that $lcm(l,m)=k''m $. Hence 
     $ {k \mid k''m} \ \Rightarrow \ {k''u^{2} \mid k''m} \ \Rightarrow \
      {u^{2} \mid m'u} \ \Rightarrow \ {u \mid m'} \ \Rightarrow \ 
      {u \mid gcd(k',m') } \ \Rightarrow \ {u=1,}$  
     contradicting the assumption that $u$ is prime. Thus $ l=uk'm' $.
\end{enumerate}
		      }
	\begin{flushright}
		$\diamondsuit$
	\end{flushright}
\end{pr}

\begin{lem}\label{cl3} %\normalsize{
	Suppose that there are $ a \in P $ and $ b \in N $ so that $ ol(a)=k$ is prime, 
	$ ol(b)=m \neq 1 $, and the following conditions hold: 
	\begin{enumerate}
	      \item $ gcd(k,m)=1 $; 
	      \item $ k \not| \; y \; ,\forall y \in olp(N) $;     
	      \item  If $ m=m'm'' $ with $gcd(m',m'')=1$ then, for each $k', \ k'' \in olp(P)$, either  
		     $$ m' \not| \; k' \quad {\rm or} \quad 
		     m'' \not| \; k''.$$ 
	\end{enumerate}
	Then these P and N do not define any circulant weighing matrix.
		 % }
\end{lem}
\begin{pr}: {\rm
	Suppose that all the conditions are satisfied. Consider the element $a-b \in \overline{\triangle} P,N $. 

	From condition (1) it follow, by Lemma \ref{cl2}, that $ ol(a-b)=km$. 
	We shall show that there is no element in 
	$ \triangle P \& \triangle N $ with orbit length equal to $ km $, thus 
	contradicting the multiset equation  
	   $$ \triangle P \& \triangle N = \overline{\triangle} P,N $$
	of claim (\ref{eq_p_n}) in Section 2.
\begin {itemize}
  \item  Suppose that $ \overline{p} \in {\triangle P} $ has $ol(\overline{p})=km$.
	 Then $ \overline{p}=p_1-p_2 $, where
	 $$ p_{1},p_{2} \in P,\; p_{1} \neq p_{2},\; ol(p_{i})=k_{i},\; i=1,2  .$$
	    By Lemma \ref{cl1}, $ km \mid lcm(k_{1},k_{2})$.
	    Since $ m \mid lcm(k_{1}, k_2)$, there exist  
			$1\leq m_1, m_2 \leq m$ 
	    such that $m=m_{1}m_{2}$, $gcd(m_1,m_2)=1$, $m_1 \mid k_1$, $m_2 \mid k_2$. 
	    This contradicts condition (3).
     \item  Suppose that $ \overline{q} \in \triangle N $ has $ol(\overline{q})=km$.
	    Then $ \overline{q}=q_{1}-q_{2} $, where
	  $$ q_{1},q_{2} \in N,\; q_{1} \neq q_{2},\; ol(q_{i})=m_{i},\; i=1,2  .$$
	    By Lemma \ref{cl1}, $ km \mid lcm(m_{1},m_{2}) $. Since $k$ is prime, 
	    This implies that either
	    $ k \mid m_{1} $ or $ k \mid m_{2} $. Contradiction with condition (2).
	      
	\begin{flushright}
		$\diamondsuit$
	\end{flushright}
\end{itemize}
	}
\end{pr}

{\rm
The above lemmas will now be used to analyze the multiset equation ((\ref{eq_p_n}) in Section \ref{sec1})
$$\triangle P \& \triangle N = \overline{\triangle} P,N .$$
A necessary condition for equality to hold is: }

{\sl For each $i$, the number of elements in $\triangle P \& \triangle N$  
with orbit length equal to $i$ is equal to the number of elements in $\overline{\triangle} P,N$ with orbit length equal to $i$.}

{\rm
For each multiset $olp(P)$, let $pol(\triangle P)$ be the set of all {\it possible orbit lengths} in  
$\triangle P$ obtained by using the above lemmas. Define similarly  $pol(\triangle N)$ from 
$olp(N)$. The following table lists the cases in which the multiset equation is false because
$$(\exists y_{0} \in \overline{\triangle} P,N ) \ 
   (\forall x \in \triangle P \& \triangle N) \ ol(x) \neq ol(y_{0}).$$ 

\large{\bf Table 2}: \normalsize{
Pairs of orbit-length partitions rejected by use of lemmas} 

\begin{flushleft}
\small{
\begin{tabular}{||l|l|l|l|l|l||}      \hline
  $\# $    &$ olp(P) $             & $pol(\triangle P) $     &$olp(N)$           &$pol(\triangle N)$       &  Rejected due to \\ \cline{1-6}
    1      & $5^{2}$               & {\scriptsize \{5\} }    & $1^{1}2^{1}3^{1}$ & {\scriptsize \{2,3,6\}} &  {\scriptsize Lemma \ref{cl2}: $k=5, \ m=2 \; {\Rightarrow} \; ol(y_{0})=10$}    \\ \hline 
    2      & $4^{1}6^{1}$          & {\scriptsize \{2,3,4,6,12\}}& $1^{1}2^{1}3^{1}$ & {\scriptsize \{2,3,6\}}                        & {\scriptsize o. k.} \\ \hline 
    3      & $3^{1}7^{1}$          & {\scriptsize \{3,7,21\}}& $1^{1}2^{1}3^{1}$ & {\scriptsize \{2,3,6\}} & {\scriptsize Lemma \ref{cl2}: $k=7, \ m=2 \; {\Rightarrow} \; ol(y_{0})=14$} \\ \hline 
    4      & $10^{1}$              & {\scriptsize \{2,5,10\}}& $1^{1}2^{1}3^{1}$ & {\scriptsize \{2,3,6\}} & {\scriptsize Lemma \ref{cl2}: $k=10, \ m=3 \; {\Rightarrow} \; ol(y_{0})=30$} \\ \hline 
    5      & $2^{1}4^{2}$          & {\scriptsize \{2,4\}}   & $3^{2}$           & {\scriptsize \{3\}}     & {\scriptsize Lemma \ref{cl2}: $k=2, \ m=3 \; {\Rightarrow} \; ol(y_{0})=6$} \\ \hline 
    6      & $1^{1}4^{1}5^{1}$     &{\scriptsize \{2,4,5,20\}}& $3^{2}$           & {\scriptsize \{3\}}     & {\scriptsize Lemma \ref{cl2}: $k=4, \ m=3 \; {\Rightarrow} \; ol(y_{0})=12$} \\ \hline 
    7      & $5^{2}$               & {\scriptsize \{5\}}      & $3^{2}$           & {\scriptsize \{3\}}     & {\scriptsize Lemma \ref{cl2}: $k=5, \ m=3 \; {\Rightarrow} \; ol(y_{0})=15$} \\ \hline 
    8      & $4^{1}6^{1}$          & {\scriptsize \{2,3,4,6,12\}}& $3^{2}$           & {\scriptsize \{3\}}                        &  {\scriptsize o. k.}    \\ \hline 
    9      & $1^{1}2^{1}7^{1}$     & {\scriptsize \{2,7,14\}} & $3^{2}$           & {\scriptsize \{3\}}     & {\scriptsize Lemma \ref{cl2}: $k=2, \ m=3 \; {\Rightarrow} \; ol(y_{0})=6$} \\ \hline 
    10     & $2^{1}8^{1}$          & {\scriptsize \{2,4,8\}}  & $3^{2}$           & {\scriptsize \{3\}}     & {\scriptsize Lemma \ref{cl2}: $k=2, \ m=3 \; {\Rightarrow} \; ol(y_{0})=6$} \\ \hline 
    11     & $1^{1}9^{1}$          & {\scriptsize \{3,9\}}    & $3^{2}$           & {\scriptsize \{3\}}                        &  {\scriptsize o. k.}    \\ \hline 
    12     & $10^{1}$              & {\scriptsize \{2,5,10\}} & $3^{2}$           & {\scriptsize \{3\}}     & {\scriptsize Lemma \ref{cl2}: $k=10, \ m=3 \; {\Rightarrow} \; ol(y_{0})=30$} \\ \hline 
    13     & $3^{2}4^{1}$          &{\scriptsize \{2,3,4,12\}}& $2^{1}4^{1}$      & {\scriptsize \{2,4\}}     & {\scriptsize Lemma \ref{cl2}: $k=3, \ m=2 \; {\Rightarrow} \; ol(y_{0})=6$} \\ \hline 
    14     & $1^{1}4^{1}5^{1}$     &{\scriptsize \{2,4,5,20\}}&$2^{1}4^{1}$       & {\scriptsize \{2,4\}}     & {\scriptsize Lemma \ref{cl2}: $k=5, \ m=2 \; {\Rightarrow} \; ol(y_{0})=10$} \\ \hline 
    15     & $5^{2}$               & {\scriptsize \{5\}}      & $2^{1}4^{1}$      & {\scriptsize \{2,4\}}     & {\scriptsize Lemma \ref{cl2}: $k=5, \ m=2 \; {\Rightarrow} \; ol(y_{0})=10$} \\ \hline 
    16     & $1^{1}3^{1}6^{1}$     & {\scriptsize \{2,3,6\}}  & $2^{1}4^{1}$      & {\scriptsize \{2,4\}}     & {\scriptsize Lemma \ref{cl2}: $k=3, \ m=4 \; {\Rightarrow} \; ol(y_{0})=12$} \\ \hline 
    17     & $4^{1}6^{1}$          & {\scriptsize \{2,3,4,6,12\}}& $2^{1}4^{1}$      &{\scriptsize \{2,4\}} &  {\scriptsize o. k.}   \\ \hline 
    18     & $3^{1}7^{1}$          & {\scriptsize \{3,7,21\}} & $2^{1}4^{1}$      & {\scriptsize \{2,4\}}     & {\scriptsize Lemma \ref{cl2}: $k=3, \ m=2 \; {\Rightarrow} \; ol(y_{0})=6$} \\ \hline 
    19     & $1^{1}9^{1}$          & {\scriptsize \{3,9\}}    & $2^{1}4^{1}$      & {\scriptsize \{2,4\}}     & {\scriptsize Lemma \ref{cl2}: $k=9, \ m=2 \; {\Rightarrow} \; ol(y_{0})=18$} \\ \hline 
    20     & $10^{1}$             &{\scriptsize \{2,5,10\}}& $2^{1}4^{1}$      &  {\scriptsize \{2,4\}}  & {\scriptsize Lemma \ref{cl6}: $k=10, \ m=4, \ u=2 \; {\Rightarrow} \; ol(y_{0})=20$} \\ \hline 
   21     &$ 3^{2}4^{1} $         & {\scriptsize \{2,3,4,12\}}& $1^{1}5^{1}$      & {\scriptsize \{5\}}     & {\scriptsize Lemma \ref{cl2}: $k=3, \ m=5 \; {\Rightarrow} \; ol(y_{0})=15$} \\ \hline 
\end{tabular}
}
\end{flushleft}

{ \ }\\

\begin{flushleft}
\small{
\begin{tabular}{||l|l|l|l|l|l||}      \hline
  $\# $    & $olp(P)$              &$pol(\triangle P) $       &$olp(N)$          & $pol(\triangle N)$      &  Rejected due to \\ \cline{1-6}
    22     & $2^{1}4^{2} $         & {\scriptsize \{2,4\}}  & $1^{1}5^{1}$      & {\scriptsize \{5\}}     & {\scriptsize Lemma \ref{cl2}: $k=2, \ m=5 \; {\Rightarrow} \; ol(y_{0})=10$} \\ \hline 
    23     & $2^{1}3^{1}5^{1}$     & {\scriptsize \{2,3,5,6,10,15\}}& $1^{1}5^{1}$      &  {\scriptsize \{5\}}& {\scriptsize o. k.} \\ \hline 
    24     &$ 5^{2}$               & {\scriptsize \{5\}}    & $1^{1}5^{1}$      & {\scriptsize \{5\}} &  {\scriptsize o. k.}\\ \hline 
    25     & $4^{1}6^{1} $         & {\scriptsize \{2,3,4,6,12\}}& $1^{1}5^{1}$ & {\scriptsize \{5\}}     & {\scriptsize Lemma \ref{cl2}: $k=4, \ m=5 \; {\Rightarrow} \; ol(y_{0})=20$} \\ \hline 
    26     & $3^{1}7^{1}$          & {\scriptsize \{3,7,21\}}& $1^{1}5^{1}$     & {\scriptsize \{5\}}     & {\scriptsize Lemma \ref{cl2}: $k=3, \ m=5 \; {\Rightarrow} \; ol(y_{0})=15$} \\ \hline 
    27     & $2^{1}8^{1}$          & {\scriptsize \{2,4,8\}}& $1^{1}5^{1}$      & {\scriptsize \{5\}}     & {\scriptsize Lemma \ref{cl2}: $k=2, \ m=5 \; {\Rightarrow} \; ol(y_{0})=10$} \\ \hline 
    28     & $10^{1}$              &{\scriptsize \{2,5,10\}}& $1^{1}5^{1}$      & {\scriptsize \{5\}} &{\scriptsize o. k.}\\ \hline 
    29     & $1^{1}2^{1}3^{1}4^{1}$&{\scriptsize \{2,3,4,6,12\}}& $6^{1}$        & {\scriptsize \{2,3,6\}} &{\scriptsize o. k.}\\ \hline    
    30     &$ 3^{2}4^{1}$          &{\scriptsize \{2,3,4,12\}}& $6^{1}$      & {\scriptsize \{2,3,6\}}  &{\scriptsize o. k.}\\ \hline 
    31     & $2^{1}4^{2}$          & {\scriptsize \{2,4\}}  & $6^{1}$      &   {\scriptsize \{2,3,6\}}& {\scriptsize Lemma \ref{cl6}: $k=4, \ m=6, \ u=2 \; {\Rightarrow} \; ol(y_{0})=12$} \\ \hline 
    32     & $2^{1}3^{1}5^{1}$     & {\scriptsize \{2,3,5,6,10,15\}}& $6^{1}$      & {\scriptsize \{2,3,6\}}     & {\scriptsize Lemma \ref{cl2}: $k=5, \ m=6 \; {\Rightarrow} \; ol(y_{0})=30$} \\ \hline  
    33     & $1^{1}4^{1}5^{1}$     & {\scriptsize \{2,4,5,20\}}& $6^{1}$      &  {\scriptsize \{2,3,6\}}     & {\scriptsize Lemma \ref{cl2}: $k=5, \ m=6 \; {\Rightarrow} \; ol(y_{0})=30$} \\ \hline 
    34     & $5^{2}$            & {\scriptsize \{5\}}& $6^{1}$      &   {\scriptsize \{2,3,6\}}     & {\scriptsize Lemma \ref{cl2}: $k=5, \ m=6 \; {\Rightarrow} \; ol(y_{0})=30$} \\ \hline 
    35     & $1^{1}3^{1}6^{1}$     &  {\scriptsize \{2,3,6\}}& $6^{1}$      &  {\scriptsize \{2,3,6\}}   & {\scriptsize o. k.}\\ \hline      
    36     & $4^{1}6^{1}$          &{\scriptsize \{2,3,4,6,12\}} & $6^{1}$      &  {\scriptsize \{2,3,6\}}  & {\scriptsize o. k.}\\ \hline
    37     & $1^{1}2^{1}7^{1}$     & {\scriptsize \{2,7,14\}}  & $6^{1}$      &  {\scriptsize \{2,3,6\}}     & {\scriptsize Lemma \ref{cl2}: $k=7, \ m=6 \; {\Rightarrow} \; ol(y_{0})=42$} \\ \hline 
    38     & $3^{1}7^{1}$          & {\scriptsize \{3,7,21\}}  & $6^{1}$      &  {\scriptsize \{2,3,6\}}     & {\scriptsize Lemma \ref{cl2}: $k=7, \ m=6 \; {\Rightarrow} \; ol(y_{0})=42$} \\ \hline 
    39     & $2^{1}8^{1}$          & {\scriptsize \{2,4,8\}}   & $6^{1}$      &{\scriptsize \{2,3,6\}} & {\scriptsize Lemma \ref{cl6}: $k=8, \ m=6, \ u=2 \; {\Rightarrow} \; ol(y_{0})=24$} \\ \hline 
    40     & $1^{1}9^{1}$          &  {\scriptsize \{3,9\}}    & $6^{1}$      &  {\scriptsize \{2,3,6\}}  & {\scriptsize Lemma \ref{cl6}: $k=9, \ m=6, \ u=3  \; {\Rightarrow} \; ol(y_{0})=18 $} \\ \hline            
    41     & $10^{1}$              & {\scriptsize  \{2,5,10\}} & $6^{1}$      &{\scriptsize \{2,3,6\}} & {\scriptsize Lemma \ref{cl6}: $k=10, \ m=6, \ u=2  \; {\Rightarrow} \; ol(y_{0}) \in \{15,30\}$} \\ \hline 

\end{tabular}
}
\end{flushleft}

\section{ Counting Arguments} \label{sec5}

%{\normalsize 
There now remain only a small number of cases.\\

\large{\bf Table 3}:{\normalsize
 \ Pairs of orbit-length partitions surviving the lemmas 
		 }
\begin{center}
{\small
\begin{tabular}{||l|l|l||}      \hline
  $\# $ &$olp(P)$               &$olp(N)$         \\ \cline{1-3}
     
     1  & $4^{1}6^{1}$          & $1^{1}2^{1}3^{1}$  \\ \hline 
     2  & $4^{1}6^{1}$          & $3^{2}$            \\ \hline 
     3  & $1^{1}9^{1}$          & $3^{2}$            \\ \hline 
     4  & $4^{1}6^{1}$          & $2^{1}4^{1}$       \\ \hline 
\end{tabular}
}
\end{center}

\begin{center}
{\small
\begin{tabular}{||l|l|l||}      \hline
  $\# $ &$olp(P)$               &$olp(N)$         \\ \cline{1-3}

     5  & $2^{1}3^{1}5^{1}$     & $1^{1}5^{1}$       \\ \hline 
     6  &$ 5^{2}$               & $1^{1}5^{1}$       \\ \hline 
     7  & $10^{1}$              & $1^{1}5^{1}$       \\ \hline 
     8  & $1^{1}2^{1}3^{1}4^{1}$ &  $6^{1}$          \\ \hline 
     9  &$ 3^{2}4^{1}$          &  $6^{1}$          \\ \hline 
     10 & $1^{1}3^{1}6^{1}$     &      $6^{1}$      \\ \hline   
     11 & $4^{1}6^{1}$          &      $6^{1}$      \\ \hline

\end{tabular}
}
\end{center}

\normalsize{
In this section we shall subject these cases to a more delicate analysis. In most cases, 
counting arguments will be used instead of simple existence considerations. 
 
Let $\langle a \rangle$ denote the $z$-orbit of $a \in {\bf Z}_n$. 
Denote also
$$ {\bf OL(\langle a \rangle - \langle a \rangle)} := \{ol(2^{i}a-2^{j}a) \mid \ i \neq j , \ 0 \leq i, \; j \leq ol(a)-1\};$$
and for $\langle a \rangle \neq \langle b \rangle$
$$ {\bf OL(\langle a \rangle - \langle b \rangle)} := \{ol(2^{i}a-2^{j}b) \mid \ 0 \leq i \leq ol(a)-1, \ 0 \leq j \leq ol(b)-1\}.$$ 

Note that $OL(\langle a \rangle - \langle a \rangle)$ and $OL( \langle a \rangle - \langle b \rangle)$ are sets 
of positive integers.  When we
write below $ol(\langle a \rangle - \langle a \rangle)$ or $ol(\langle a \rangle - \langle b \rangle)$ we mean an arbitrary element of the  
corresponding set. We shall now analyze all the cases in the above table, one by one.

\begin{enumerate}
  \item  $ P=\{a,2a,4a,8a,b,2b,4b,8b,16b,32b\}; \  N=\{c,d,2d,e,2e,4e\} $; \\  
	 $ ol(a)=4, \ ol(b)=6; \  ol(c)=1, \  ol(d)=2,\  ol(e)=3 $.

	 Let us count the elements with orbit length 12 in $ \triangle P \& \triangle N$ and 
	 in $ \overline{\triangle}P,N$.               
	 By Lemma \ref{cl6} 
	   $ ol(\langle a \rangle - \langle b \rangle)=12 $ and by Lemma \ref{cl2} $ ol(\langle a \rangle - \langle e \rangle)=12 $. 
	 It easy to see that no other combinations of the orbit lengths in this case yields 12. Therefore the number of elements with orbit length 12 in $\triangle P \& \triangle N$ is 
 $2 \times 4 \times 6=48$ and in $\overline{\triangle}P,N$ it is $2 \times 4 \times 3=24$.
	  Thus, we obtain the contradiction
	$$ \overline{\triangle}P,N \neq \triangle P \& \triangle N.$$
  
  \item  $ P=\{a,2a,4a,8a,b,2b,4b,8b,16b,32b\}; \  N=\{c,2c,4c,d,2d,4d\} $; \\  
	 $ ol(a)=4, \ ol(b)=6; \  ol(c)=3, \  ol(d)=3 $.
	 $$\triangle P \& \triangle N:$$
	 $$\triangle P : \ 
	  ol(\langle b \rangle - \langle b \rangle) \in \{2,3,6\}, \mbox{ by Lemma \ref{cl1} and Observation \ref{cl01}}.$$
	 $$\triangle N :  \  
	 ol(\langle c \rangle - \langle c \rangle)=ol(\langle c \rangle - \langle d \rangle)=ol(\langle d \rangle - \langle d \rangle)=3, \mbox{ by Corollary \ref{cl11}}.$$
	 $ \overline{\triangle} P,N$: 
	 $$ ol(\langle a \rangle - \langle c \rangle)=ol(\langle a \rangle - \langle d \rangle)=12,  \mbox{ by Lemma \ref{cl2}} .$$
	 $$ ol(\langle b \rangle - \langle c \rangle), \ ol(\langle b \rangle - \langle d \rangle) \in \{2,6\},  \mbox{ by Lemma \ref{cl4}}.$$
	 Thus there is no element in $ \overline{\triangle} P,N$  whose orbit length is equal to 3.
	 Hence, we conclude that $ \overline{\triangle} P,N  \neq \triangle P \& \triangle N $.\\

  \item  $ P=\{a,b,2b,4b,8b,16b,32b,64b,128b,256b\}; \  N=\{c,2c,4c,d,2d,4d\} $; \\  
	 $ ol(a)=1, \ ol(b)=9; \  ol(c)=3, \  ol(d)=3 $.
	 
	 $$\triangle N : \ 
	  ol(\langle c \rangle - \langle c \rangle)=ol(\langle c \rangle - \langle d \rangle )=ol(\langle d \rangle - \langle d \rangle)=3,  \mbox{ by Corollary \ref{cl11}}$$
	 and the number of these elements is $3\times2+2\times3\times3+3\times2=30$;
 
	 $$ \overline{\triangle} P,N :  \ 
	  ol(\langle a \rangle - \langle c \rangle)=ol(\langle a \rangle - \langle d \rangle)=3, \mbox{ by Observation \ref{cl0}}$$   
	 and the number of these elements is $2\times1\times3+2\times1\times3=12$;
	 $$ ol(\langle b \rangle - \langle c \rangle) = ol(\langle b \rangle - \langle d \rangle)=9, \mbox{ by Lemma \ref{cl4}} .$$
	 Thus the number of elements in $ \overline{\triangle} P,N $ whose orbit length is equal to 3 
	 is smaller then the number of such elements in $ \triangle P  \& \triangle N $.
	 Hence, we obtain that $ \overline{\triangle} P,N  \neq \triangle P  \& \triangle N  $.

  \item  $ P=\{a,2a,4a,8a,b,2b,4b,8b,16b,32b\}; \  N=\{c,2c,d,2d,4d,8d\} $; \\  
	 $ ol(a)=4, \ ol(b)=6; \  ol(c)=2, \  ol(d)=4 $.\\
	 In this case we do not get a contradiction by using the lemmas from the 
	 Section~\ref{sec4}, and this case remains as a candidate to be dealt with in the 
	 next section.

  \item  $ P=\{a,2a,b,4b,8b,c,2c,4c,8c,16c\}; \  N=\{d,e,2e,4e,8e,16e\} $; \\  
	 $ ol(a)=2, \ ol(b)=3, \  ol(c)=5; \  ol(d)=1, \  ol(e)=5 $.\\
	 $$\triangle P : \ 
	  ol(\langle a \rangle - \langle b \rangle)=6, \mbox{ by Lemma \ref{cl2}}. $$
	 It is easy to check that 
	 %$ \overline{\triangle} P,N $: 
	 %$$ ol(a-d)=ol(d-a)=2, \mbox{ by Observation \ref{cl0}};$$
	 %$$ ol(b-d)=ol(d-b)=3, \mbox{ by Observation \ref{cl0}};$$
	 %$$ ol(c-d)=ol(d-c)=5, \mbox{ by Observation \ref{cl0}};$$
	 %$$ ol(a-e)=ol(e-a)=10, \mbox{ by Lemma \ref{cl2}};$$
	 %$$ ol(b-e)=ol(e-b)=15, \mbox{ by Lemma \ref{cl2}};$$
	 %$$ ol(c-e)=ol(e-c)=5, \mbox{ by Observation \ref{cl01}   };$$
	 there is no element in $ \overline{\triangle} P,N $ whose orbit length is equal to 6.
	 Hence, we conclude that $ \overline{\triangle} P,N  \neq \triangle P \& \triangle N $.
  
  \item  $ P=\{a,2a,4a,8a,16a,b,2b,4b,8b,16b\}; \  N=\{c,d,2d,4d,8d,16d\} $; \\ 
	 $ ol(a)=5, \ ol(b)=5; \  ol(c)=1, \  ol(d)=5 $.\\
	 None of the lemmas leads to a contradiction in this case.
	 Therefore this case still remains as a candidate.

  \item  $ P=\{a,2a,4a,8a,16a,32a,64a,128a,256a,511a\}; \  N=\{b,c,2c,4c,8c,16c\} $; \\  
	 $ ol(a)=10; \ ol(b)=1, \  ol(c)=5 $.\\  
	 $$\triangle N: \ 
	 ol(\langle b \rangle - \langle c \rangle)=5, \mbox{ by Observation \ref{cl0}}.$$
	 $ \overline{\triangle} P,N $: 
	 $$ ol(\langle a \rangle - \langle b \rangle)=10, \mbox{ by Observation \ref{cl0}} ;$$
	 $$ ol(\langle a \rangle - \langle c \rangle) \in \{2,10\}, \mbox{ by Lemma \ref{cl4}}.$$
	 Thus there is no element in $ \overline{\triangle} P,N $ with orbit length equal to 5.
	 Hence, we conclude that $ \overline{\triangle} P,N  \neq \triangle P \& \triangle N $.
 
  \item  $ P=\{a,b,2b,c,2c,4c,d,2d,4d,8d\}; \  N=\{e,2e,4e,8e,16e,32e\}; $ \\  
	 $ ol(a)=1, \ ol(b)=2, \  ol(c)=3, \  ol(d)=4; \  ol(e)=6 $.\\
	 $$\triangle P: \ ol(\langle c \rangle - \langle d \rangle)=12, \mbox{ by Lemma \ref{cl2}},$$  
	 and these are the only elements in $\triangle P$ with orbit length 12. 
	 Clearly, there are no elements of orbit length 12 in $\triangle N$. 
	 $$ \overline{\triangle} P,N: \ 
	  ol(\langle d \rangle - \langle e \rangle)=12 , \mbox{ by Lemma \ref{cl6}},$$
	 and there are no other                  
	 such elements in $\overline{\triangle}P,N$. Thus, we obtain that 
	 in $ \overline{\triangle} P,N $ there are $ 2 \times (4 \times 6)=48 $ elements
	 whose orbit length is equal to 12, while in $\triangle P \& \triangle N $
	 there are only $2 \times (3 \times 4)=24$ such elements. Hence, \- 
	 $ \overline{\triangle} P,N \neq \triangle P \& \triangle N $.
  
  \item  $ P=\{a,2a,4a,b,2b,4b,c,2c,4c,8c\}; \  N=\{d,2d,4d,8d,16d,32d\} $; \\  
	 $ ol(a)=3, \ ol(b)=3, \  ol(c)=4; \  ol(d)=6 $.\\
	 $$\triangle P : \ ol(\langle a \rangle - \langle a \rangle)=ol(\langle a \rangle - \langle b \rangle)=ol(\langle b \rangle - \langle b \rangle)=3, \mbox{ by Corollary \ref{cl11}},$$
	 %$ \overline{\triangle} P,N$:
	 %$$ ol(d-a)=ol(a-d)=ol(d-b)=ol(b-d) \in \{2,6\}, \mbox{ by Lemma \ref{cl4}}; $$
	 %$$ ol(d-c)=ol(c-d)=12 , \mbox{ by Lemma \ref{cl6}};$$ 
	 but there is no element in $ \overline{\triangle} P,N$ whose orbit length is
	 equal to 3 (may be shown using lemmas \ref{cl4} and \ref{cl6}). Hence, $ \overline{\triangle} P,N \neq \triangle P \& \triangle N $.

  \item  $ P=\{a,b,2b,4b,c,2c,4c,8c,16c,32c\}; \  N=\{d,2d,4d,8d,16d,32d\} $; \\  
	 $ ol(a)=1, \ ol(b)=3, \  ol(c)=6; \  ol(d)=6 $.\\
	 None of the lemmas leads to a contradiction in this case, and it remains as 
	 a candidate.

  \item  $ P=\{a,2a,4a,8a,b,2b,4b,8b,16b,32b\}; \  N=\{c,2c,4c,8c,16c,32c\} $; \\  
	 $ ol(a)=4, \ ol(b)=6; \  ol(c)=6 $.\\
	 $$\triangle P: \ 
	 ol(\langle a \rangle - \langle a \rangle) \in \{2,4\}, \mbox{ by Lemma \ref{cl1} and Observation \ref{cl01}}.$$
	 %$ \overline{\triangle} P,N$: 
	 %$$ ol(a-c)=ol(c-a)=12, \mbox{ by Lemma \ref{cl6}} ;$$
	 %$$ ol(b-c)=ol(c-b) \in \{2,3,6\}, \mbox{ by Lemma \ref{cl1} and Observation \ref{cl01}} ;$$
	 Note that $ol(2a-a)=ol(a)=4$, so there is at least one element of orbit length
	 4 in $\triangle P $ and therefore in $\triangle P \& \triangle N $. On the other hand,  
	 there are no elements of orbit length 4 in $ \overline{\triangle} P,N$.
	 Hence, $ \overline{\triangle} P,N \neq \triangle P \& \triangle N $.

\end{enumerate}

%}

\section{Final Analysis} \label{sec6}

%\normalsize{
In this section we shall analyze more closely the three cases that survived 
the previous inspection (cases 4, 6 and 10). Let us list them again, in a 
different order: 
\begin{enumerate}
      \item $olp(P)=5^{2}, \ olp(N)=1^15^1$.
      \item $olp(P)=1^13^16^1, \ olp(N)=6^1$.
      \item $olp(P)=4^16^1, \ olp(N)=2^14^1$. 
\end{enumerate}

Before embarking upon the detailed, examination of these cases, we need some 
general theorems. 

\begin{thm} \label{tm1}
%{\normalsize    
 Suppose that $A$ is a $v \times v$ weighing matrix with $v=km$ 
     which has the block form: 
	  \[ \left( \begin{array}{cccc}
	W   & 0 & \ldots & 0 \\
	0 & W & \ldots & 0 \\
	\vdots  \\
	0 & 0 & \ldots & W \\
	       \end{array}      \right)   \]
	  where $W$ is a circulant $m \times m$ matrix. 
	  Then there is a $v \times v$ permutation matrix $P$ such that 
	  $P^{-1}AP$ is a circulant weighing matrix.
%}
\end{thm}     
\begin{pr}:
{\rm
	\ Let $A=(a_{ij})$ and $W=(w_{ij})$ be the above matrices. 
	Let $P=(p_{ij})_{i,j=0}^{v-1}$ be the $v \times v$ permutation 
	matrix defined by: 
$$ p_{ij}=1 \quad \Longleftrightarrow \quad i=rm+s, \ j=sk+r  \quad 
	{\rm for \ some} \quad 0 \leq r \leq k-1, \quad 0 \leq s \leq m-1.$$
	Multiplying the matrix $A$ by $P^{-1}$ on the left turns row  
	$ rm+s $ of the matrix $A$ into row $sk+r$, while 
	multiplying $A$ by $P$ on the right turns 
	column  $ rm+s $ of the matrix $A$ into column $sk+r$.    
	We'll prove now that $B:=P^{-1}AP$ is a circulant weighing matrix.  
	Because permutation of rows and columns of a weighing matrix gives a weighing 
	matrix, we only have to prove that $B$ is circulant. 
	In order to do so, we will prove that  
		$$ b_{ij}=b_{i+1,j+1} \quad (\forall \quad  0 \leq i,j \leq v-1),$$ 
	where addition of indices is modulo $v$.
	Let 
		$$ i=s_1k+r_1 \quad {\rm and } \quad j=s_2k+r_2 \quad
		 (0 \leq r_{1},r_{2} \leq k-1, \quad 0 \leq s_{1},s_{2} \leq m-1).$$
	Then 
		$$ i+1=s_1k+r_1+1 \quad {\rm and } \quad j+1=s_2k+r_2+1 .$$
  The following table shows which rows and columns of the matrix $A$ 
  correspond to given rows and columns of the matrix $B$.\\
  { \ }\\
\begin{tabular}{||l|lr||}      \hline
  Matrix B     &  Matrix A & \\ \cline{1-3}
  row $i$      &  row  $r_1m+s_1$   & \\ \hline
  column $j$   & column  $r_2m+s_2$ & \\ \hline
  row $i+1$    & row  $(r_1+1)m+s_1$,        &  if $r_1 \neq k-1$    \\ \cline{2-3}
	       & row  $0 \cdot m+(s_1+1)$,   &  if $r_1=k-1, \quad s_1 \neq m-1$         \\ \cline{2-3}
	       & row  $0 \cdot m+0$,   &  if $r_1=k-1, \quad s_1=m-1$         \\ \hline
  column $j+1$ & column  $(r_2+1)m+s_2$,     &  if $r_2 \neq k-1$    \\ \cline{2-3}
	       & column  $0 \cdot m+(s_2+1)$,&  if $r_2=k-1, \quad s_2 \neq m-1$          \\ \cline{2-3}
	       & column  $0 \cdot m+0$,&  if $r_2=k-1, \quad s_2=m-1$          \\ \hline
\end{tabular}
{ \ }\\
{ \ }\\
     
     Thus the following cases are possible:
	
\begin{itemize}
  \item $\underline {r_1=r_2}$.\\
	In this case row $r_1m+s_1$ of $A$ is row $s_1$ in diagonal block number 
	$r_1$, and column $r_2m+s_2$ is column $s_2$ in the same block.  
	Hence
		$$b_{ij}=a_{(r_1m+s_1),(r_2m+s_2)}=w_{s_1s_2}.$$
	If $r_1 \neq k-1$ then
		$$b_{(i+1),(j+1)}=a_{((r_1+1)m+s_1),((r_2+1)m+s_2)}=w_{s_1s_2}.$$           
	Otherwise,  $r_1=k-1$ and 
	   $$b_{(i+1),(j+1)}=a_{(0 \cdot m+(s_1+1)),(0 \cdot m +(s_2+1))}=w_{(s_1+1),(s_2+1)}=w_{s_1s_2}.$$           
	The last equality follows from $W$ being circulant.
	Here $s_1+1$ and $s_2+1$ are taken modulo $m$, covering also the 
	cases  where 
	  $$s_1=m-1 \quad {\rm or} \quad s_2=m-1 \quad ({\rm or \ both}).$$
	Hence in all cases $b_{ij}=b_{(i+1)(j+1)}$.         

  \item $\underline {r_1 \neq r_2}$.\\
	In this case the entry in row $r_1m+s_1$ and column $r_2m+s_2$ 
	does not belong to a diagonal block of $A$. 
	
	%(and also the row $(r_1+1)m+s_1$ and the column $(r_2+1)m+s_2$;  
	%the row $(r_1+1)m+s_1$ and the column $s_2+1$) 
	%are the row and the column in the distinct blocks.  
	
	Hence
		$$b_{ij}=a_{(r_1m+s_1),(r_2m+s_2)}=0.$$
	If $r_1 \neq k-1$ and $r_2 \neq k-1$ then, similarly, 
		$$b_{(i+1),(j+1)}=a_{((r_1+1)m+s_1),((r_2+1)m+s_2)}=0.$$           
	Otherwise, with no loss of generality suppose that 
	$r_1 \neq k-1$ and $r_2=k-1$. then
	$$b_{(i+1),(j+1)}=a_{((r_1+1)m+s_1),(0 \cdot m +(s_2+1))}=0.$$           
	Hence again $b_{ij}=b_{(i+1),(j+1)}$.         
\end{itemize}

 We have proved that
     $$ b_{ij}=b_{i+1,j+1} \quad (\forall \quad  0 \leq i,j \leq v-1)$$ 
 and therefore the weighing matrix $B$ is circulant.  
		   \begin{flushright}
				$\diamondsuit$
			  \end{flushright}
}
\end{pr}

\begin{thm} \label{theor1}
%      {\normalsize
       If $CW(n,k) \neq \emptyset$ then $CW(mn,k) \neq \emptyset$ for every $m \geq 1$.
%         }
\end{thm}
\begin{pr}:{\rm
	\ Let $W \in CW(n,k)$ and let $I_m$ be the identity matrix of order $m \geq 1$.
	Note that $I_m$ is a circulant weighing matrix of order $m$ and 
	weight 1.
	Hence the Kronecker product of these matrices gives the matrix
	$W'=I_m \otimes W \in W(mn,k)$: 
	\[ W'= \left( \begin{array}{cccc}
	W  & 0 & \ldots & 0  \\
	0  & W & \ldots & 0  \\            
	\vdots  \\
	0  & 0 & \ldots & W  \\
	       \end{array}      \right)   \]

	By the previous Theorem \ref{tm1} there is an $mn \times mn$ 
	permutation matrix $P$ such that $P^{-1}W'P \in CW(mn,k)$. 
	Thus $CW(mn,k) \neq \emptyset$.
		   
		   \begin{flushright}
				$\diamondsuit$
		   \end{flushright}
     
     }
\end{pr}

In the sequel we shall attempt to prove a converse to Theorem 7.2, but this 
will be done separately for each of the cases at hand. In each case we shall 
assume a specific pairs $(olp(P),olp(N))$. 

%From this theorem we get the following conclusion.

%\begin{con} \label{con1}
%%  {\normalsize  
%    It is enough to check  
%    the existence of circulant weighing matrices (with given weight $k$ and 
%    orbit length partitions) for some ``minimal" order $n$. 
%    The result (if positive) will also imply the 
%    existence of circulant weighing matrices of order $mn$ 
%    and weight k (and the same orbit length partitions), for all $m \geq 1 $.
%%  }
%\end{con}

We now proceed with the analysis of the above three cases.
\begin{enumerate}
  \item  $ P=\{a,2a,4a,8a,16a,b,2b,4b,8b,16b\}; \  N=\{c,d,2d,4d,8d,16d\} $; \\ 
	 $ ol(a)=5, \ ol(b)=5; \  ol(c)=1, \  ol(d)=5 $.\\
	 
	 $P$ and $N$ contain orbits of lengths 1 and  5.
	 In Section \ref{sec3} we found necessary and sufficient conditions on 
	 $n$ for the existence of an element in ${\bf Z}_n$  
	 with orbit length equal to $i$, for each $1 \leq i \leq 6$.  
	 In the present case, $n$ must satisfy the conditions  for $i=1$ and 
	 $i=5$.  
	 \begin{itemize}
	     \item \underline{$i=1$}: $n$ arbitrary.
	     \item \underline{$i=5$}: $n$ must be divisible by 31.
	 \end{itemize}    
	 We may thus assume that $n=31m$ for some (odd) integer $m$. 
	 We shall now state and prove a converse to Theorem 7.2, especially for 
	 the current case.  

\begin{thm} \label{th1}
 %     {\normalsize
{ \ }\\
(i)    For each odd $m \geq1 $, if $w(x) \in CW(31m,16)$ has (for the 
       multiplier $t=2$)
		$$olp(P)=5^2, \quad olp(N)=1^15^1$$ 
	then there exists $w_0(x) \in CW(31,16)$ 
       such that 
		$$ w(x)=w_0(x^m).$$
	   
(ii)    If $w_0(x^m), \ \widetilde{w}_0(x^m)$ are equivalent in $CW(31m,16)$  
	then $w_0(x), \ \widetilde{w}_0(x)$ are equivalent in $CW(31,16)$  
	   %if type($CW(31,16)$)=$m \times $type($CW(31m,16)$) then  
       %$$CW(31,16) \neq \emptyset \ {\rm iff} \ CW(31m,16) \neq \emptyset.$$
  %       }
\end{thm}
\begin{pr}:{\rm
{ \ }\\        

(i)     Let $m \geq 1$ be an odd integer and assume that $w(x) \in CW(31m,16)$ 
	with the given $olp(P), \ olp(N)$. 
	Then the unique element with orbit length equal to 1 is $0 \in N$. Let 
	$x \in N$ and $y, \; z \in P$ be generators for the three orbits of length 5. 
	According to Section \ref{sec3},
	  $$ x=mk $$
	for some $1 \leq k \leq 30$. Similarly, 
	       $$ y=mk', \quad z=mk'' $$
	for some $1 \leq k', \ k'' \leq 30$. 
Thus 
	$$ m | s \quad (\forall s \in P \cup N)$$
	 so that there is a (unique) polynomial $w_0(x) \in {R}_{31}$ 
	 such that 
	 $$ w(x)=w_0(x^{m}).$$
	 Clearly, $w_0(x) \in CW(31,16)$.

(ii)    Let $w_0(x), \ \widetilde{w}_0(x) \in CW(31,16)$ be such that the 
	polynomials $w_0(x^m), \ \widetilde{w}_0(x^m) \in CW(31m,16)$ are equivalent. 
	Thus there exist  $s \in {\bf Z}_{31m}$ and $t \in {\bf Z}^{*}_{31m}$ 
	such that 
		$$ \widetilde{w}_{0}(x^m)=x^{s}w_0(x^{mt}) \quad ({\rm 
in} \ R_{31m}).$$ 
	All the powers of $x$ with non-zero coefficients in $\widetilde{w}_{0}(x^m)$ 
	or $\widetilde{w}_{0}(x^{mt})$ are divisible by $m$. Therefore $s=ms_1$ 
	for a suitable $s_1 \in {\bf Z}_{31}$, and we conclude that 
	     $$ \widetilde{w}_{0}(x)=x^{s_1}w_0(x^{t}) \quad ({\rm \ in \ } R_{31}).$$ 
	Note that $t \in {\bf Z}^{*}_{31m}$ may also be viewed as $t \in {\bf Z}^{*}_{31}$. 
		   \begin{flushright}
				$\diamondsuit$
		   \end{flushright}
     
     }
\end{pr}

	 We now face the problem of finding all $w(x) \in CW(31,16)$ 
	 with the above $olp(P)$ and $olp(N)$, and sorting them into equivalence 
	 classes. The data that we have are
	      $$n=31, \quad k=16, \quad t=2, \quad olp(P)=5^2, \quad olp(N)=1^{1}5^1.$$ 
	 We will find $P$ and $N$ with the help of a computer program.               
	 We are looking for 
	 $w(x)=w_{0}+w_{1}x+\cdots+w_{30}x^{30} \in CW(31,16)$ with $w_{0}=-1$ 
	 ($0 \in N$ since only 0 has orbit length equal to 1). The indices $i$ for all 
	 the other nonzero $w_i$ belong to orbits of length 5. 
	 According to Section \ref{sec3}, there are six different orbits of length 5
	 in ${\bf Z}_{31}$. For our $w(x)$ we need three of them. 
	 Thus the Pascal program must check $6{5 \choose 2}=60$ cases.
	 Each case gives explicit $P$ and $N$ and therefore also $w(x) \in R_{31}$ 
	 which defines a circulant $\{0,1,-1\}$-matrix. 
	 In order to verify that this is a weighing matrix, we have to check 
	 the following conditions:
	 \begin{itemize}
	\item $\sum_{i=0}^{30}w^2_i=16 $;
	\item $\sum_{i=0}^{30}w_{i}w_{i+j} =0, \quad  1 \leq j \leq 30$;
	 \end{itemize}
	 Note that the first condition is automatically satisfied (since 
	 $w_i=\pm 1$); and in the second condition it is sufficient to check only cases up to 
	 $j=\lfloor \frac{31}{2} \rfloor$, since the other values of $j$ are complements of 
	 these to 31. Hence for every candidate $w(x)$
	 the program should make 15 comparisons. 
	 
	 This program was written and run, finding altogether 12 solutions. 
	The list of solutions in $CW(31,16)$ includes all those  
	previously obtained by R. Eades \cite{ead80} and by Y. Strassler 
	\cite{str95}. Another program designed to check possible equivalence 
	between the matrices obtained. It showed that every $w(x) \in CW(31,16)$ 
	is equivalent to one of the following two: 
	 \begin{enumerate}
	      \item $w_{1}(x)=-1-x^{1}-x^{2}+x^{3}-x^{4}+x^{6}+x^{7}-x^{8}+x^{12}+x^{14}-x^{16}+x^{17}+x^{19}+x^{24}+x^{25}+x^{28}, $
	      \item $w_{2}(x)=-1-x^{1}-x^{2}-x^{4}+x^{5}-x^{8}+x^{9}+x^{10}+x^{15}-x^{16}+x^{18}+x^{20}+x^{23}+x^{27}+x^{29}+x^{30}. $
	 \end{enumerate}
	 
	It is easy to see that $w_{1}(x)$ and $w_{2}(x)$ are inequivalent. 
	Indeed, by Section \ref{sec3} there are exactly 6 different orbits of 
	length 5 in $ {\bf Z}_{31} $. These are 
	$$C_0=\{1,2,4,8,16\},$$ 
	$$C_1=\{3,6,12,24,17\},$$
	$$C_2=\{5,10,20,9,18\},$$
	$$C_3=\{7,14,28,25,19\},$$
	$$C_4=\{11,22,13,26,21\},$$
	$$C_5=\{15,30,29,27,23\}.$$
	Denote by $C_{\infty}=\{0\}$ the unique orbit of length 1.

	Thus we obtain
\begin{itemize}
 \item $w_{1}(x)$:
  $$ P_1=\{ 3,6,7,12,14,17,19,24,25,28\}= C_1 \cup C_3,$$
  $$ N_1=\{ 0,1,2,4,8,16\}=C_{\infty} \cup C_0;$$

 \item $w_{2}(x)$:
  $$ P_2=\{5,9,10,15,18,20,23,27,29,30\}= C_2 \cup C_5,$$
  $$ N_2=\{ 0,1,2,4,8,16\}=C_{\infty} \cup C_0.$$
\end{itemize}    
			 
	 Assume that $w_{1}(x)$ and $w_{2}(x)$ are equivalent:  
	 $$ w_2(x)=x^{s}w_1(x^t)  \qquad (s \in {\bf Z}_{31}, \ t \in {\bf Z}^{*}_{31}).$$ 
	 Since 2 is a fixing multiplier for both $w_1(x)$ and $w_2(x)$, 
		$$ w_1(x^2)=w_1(x) \ {\rm and } \ w_2(x^2)=w_2(x).$$ 
	 Thus
	   $$w_2(x)=w_2(x^2)=x^{2s}w_1(x^{2t})=x^{2s}w_1(x^t)=x^{s}w_2(x^t).$$ 
	 It follows that 
		$$ P_2=s+P_2 \ {\rm and } \ N_2=s+N_2,$$
	 and obviously this implies $s=0$ (for the given $P_2$ and $N_2$). 
	 Thus  
		$$ w_2(x)=w_1(x^t)$$                 
	 so that 
		$$ P_2=tP_1 \ {\rm and } \ N_2=tN_1.$$
	 Multiplication by $t$ maps 2-orbits to 2-orbits (of the same length, 
	 and thus $tC_{\infty}=C_{\infty}$ and $tC_0=C_0$. It follows that 
	 $t=t \cdot 1 \in tC_{\infty}=C_{\infty}$, i.e., $t$ is a power of 
	 $2 (mod 31)$. Thus $w_2(x)=w_1(x^t)=w_1(x)$. 

	 By Theorems \ref{theor1} and \ref{th1} we get that for every odd  
	 $d \geq 1$ there are 
	 two distinct equivalence classes in $CW(31m,16)$ with $olp(P)=5^2$ 
	 and $olp(N)=1^15^1$. They are
	 \begin{itemize}   
	   \item
	   $ \widetilde{w_1}(x)=-1-x^{m}-x^{2m}+x^{3m}-x^{4m}+x^{6m}+x^{7m}-x^{8m}+x^{12m}+x^{14m}-x^{16m}+x^{17m}+x^{19m}+x^{24m}+x^{25m}+x^{28m}, $
	   \item
	   $ \widetilde{w_2}(x)=-1-x^{m}-x^{2m}-x^{4m}+x^{5m}-x^{8m}+x^{9m}+x^{10m}+x^{15m}-x^{16m}+x^{18m}+x^{20m}+x^{23m}+x^{27m}+x^{29m}+x^{30m}. $
	 \end{itemize}

{ \ }\\

  \item  $ P=\{a,b,2b,4b,c,2c,4c,8c,16c,32c\}; \  N=\{d,2d,4d,8d,16d,32d\} $; \\  
	 $ ol(a)=1, \ ol(b)=3, \  ol(c)=6; \  ol(d)=6 $.\\
	 
	 $P$ and $N$ contain orbits of lengths 1, 3 and 6.
	 By Section \ref{sec3}, $n$ must satisfy the following conditions  for 
	 the existence of orbits of length $i$, for $i \in \{1,3,6\}$.
	 By Section \ref{sec3}
	 \begin{itemize}
	     \item \underline{$i=1$}: $n$ arbitrary.
	     \item \underline{$i=3$}: $n$ must be divisible by 7.
	     \item \underline{$i=6$}: In this case, the precise restrictions on $n$ depend on the number
	     of different orbits of length 6 in $P \cup N$. We need two   
	     orbits. Hence there are two possibilities: 
	     \begin{enumerate}
	\item $63 \mid n$, so there are 9 different orbits of length 6.
	     Therefore there are $9 \times 8=72$ possibilities for the choice of the 
	     two orbits of length 6 in $P \cup N$. 
	\item $21 \mid n$ but $63 \not| \; n$, and then there are exactly 
	     two different orbits of length 6.
	     Therefore there are two possibilities for the choice of the two 
	     (ordered) orbits of length 6 in $P \cup N$.
	     
	     %Note that this case includes in itself the first one (a).
	     \end{enumerate}
	 
	 \end{itemize}
	 (Recall that if $9 \; | \; n$ but $63\not| \; n$ then there is only one orbit 
	 of length 6.) 
	 Hence in the present case necessarily $21 \mid n$.

\begin{thm} \label{th2}
      %{\normalsize
{ \ }\\
(i)    For each odd $m \geq1 $, if $w(x) \in CW(21m,16)$ has (for the 
       multiplier $t=2$) 
		$$olp(P)=1^13^16^1, \quad olp(N)=6^1$$ 
       then:
       \begin{itemize}
	\item If $3 \not| \; m$ then there exists 
		$$w_0(x) \in CW(21,16) \quad {\rm s. \ t.} \quad w(x)=w_0(x^m).$$ 
	\item If $3 | m$ then there exists  
	$$w'_0(x) \in CW(63,16) \quad {\rm s. \ t.} \quad w(x)=w'_0(x^{\frac{m}{3}}).$$ 
	\end{itemize}
	
(ii)  If $w(x)$ and $\widetilde{w}(x)$ are equivalent in $R_{21m}$ then: 
	\begin{itemize}
	  \item 
		If  $3 \not| \; m$ then $w_0(x)$ and $\widetilde{w}_0(x)$ are 
		equivalent in $R_{21}$. 
	  \item 
		If $3 | m$ then $w'_0(x)$ and $\widetilde{w}'_0(x)$ are 
		equivalent in $R_{63}$. 
		
	\end{itemize}
       % }
\end{thm}
\begin{pr}
{\rm
{ \ }\\

(i)     Let $m \geq 1$ be an odd integer and assume that $w(x) \in CW(21m,16)$ with the above $olp(P), \ olp(N)$.
	Then $0 \in P$ (the unique element with orbit length equal to 1).  
	Let $ x \in P$ be a generator of the orbit of length 3. 
	According to Section \ref{sec3} 
	  $$ \exists  j \in \{1,\cdots,6\} : \quad x=\frac{21mj}{7}=3mj \quad \Rightarrow \quad m | x.$$
	Let $y \in P$ and $z \in N$ be generators of the orbits of length 6. 
	Therefore by the results of Section \ref{sec3}
	  $$ \exists k \in \{1, \cdots, 62\}, \quad 9 \not| \; k \ {\rm and} \ 21 \not| \; k $$
	such that
			$$ y=\frac{21mk}{63}=\frac{mk}{3}.$$       
       Similarly for $z$. 
       The following cases are possible: 
       \begin{itemize}
	   \item $3 \not| \; m$.\\ 
	     %Because $CW(21m,16) \neq \emptyset$ then must be 
	     Here necessarily $3 | k$.                  
	     Hence $m | y$ and similarly for $z$. Thus we obtain 
		$$m|s \quad (\forall s \in P \cup N),$$ 
	     so
	     there is a (unique) polynomial $w_0(x) \in {R}_{21}$ s.t. 
		$w(x)=w_0(x^m)$.
	     Obviously, $w_0(x) \in CW(21,16)$.
	   \item $3 | m$.\\
	     Let $m=3m'$. Then $ y=m'k \quad \Rightarrow \quad m' | y.$   
	     Similarly for $z$. Note that also $m' | x$. Thus 
	      $$m'|s \quad  (\forall s \in P \cup N),$$
	     so there is a (unique) polynomial $w'_0(x) \in {R}_{63}$ s.t. 
		$$ w(x)=w'_0(x^{m'})=w'_0(x^{\frac{m}{3}}).$$
	     Obviously, $w'_0(x) \in CW(63,16)$.

	\end{itemize}

(ii)    Proof similar to that of Theorem \ref{th1}(ii). 
		   \begin{flushright}
				$\diamondsuit$
		   \end{flushright}
}
\end{pr}     
     
	 By Theorem \ref{theor1} and Theorem \ref{th2} 
	 we now need to find all equivalence classes in $CW(21,16)$ and 
	 in $CW(63,16)$ (with the given $olp(P)$ and $olp(N)$). Note that if 
	 $w(x) \in CW(21,16)$ then $w(x^3) \in CW(63,16)$.  
	 The data that we have consist of 
		$$n=63, \quad k=16, \quad t=2, \quad olp(P)=1^13^16^1, \quad olp(N)=6^1.$$ 
	 We will search for $P$ and $N$ with the help of a Pascal program.               
	 This program is very similar to the one 
	 described above in the case of $CW(31,16)$. It was run, giving 8 
	 solutions. Another Pascal program was designed to check equivalence 
	 between the polynomials obtained. It showed that every $w(x) \in CW(63,16)$  
	 is equivalent to one of the following two polynomials: 
	 \begin{enumerate}
	\item $w_{1}(x)=1-x^{1}-x^{2}-x^{4}-x^{8}+x^{9}+x^{13}-x^{16}+x^{18}+x^{19}+x^{26}-x^{32}+x^{36}+x^{38}+x^{41}+x^{52} $
	
	\item $w_{2}(x)=1-x^{3}-x^{6}-x^{12}+x^{15}-x^{24}+x^{27}+x^{30}-x^{33}+x^{39}+x^{45}-x^{48}+x^{51}+x^{54}+x^{57}+x^{60}  $
	 \end{enumerate}
	 
	As in the case of $CW(31,16)$, it may be easily shown 
	that $w_1(x)$ and $w_2(x)$ are inequivalent. 

	By Theorems \ref{theor1} and \ref{th2} 
	 there are two distinct equivalence classes in $CW(63m,31)$, for each 
	 odd $m \geq 1$: 
	 \begin{itemize}   
	 \item
	   $ \widetilde{w_1}(x)=1-x^{m}-x^{2m}-x^{4m}-x^{8m}+x^{9m}+x^{13m}-x^{16m}+x^{18m}+x^{19m}+x^{26m}-x^{32m}+x^{36m}+x^{38m}+x^{41m}+x^{52m}, $
	\item 
	   $ \widetilde{w_2}(x)=1-x^{3m}-x^{6m}-x^{12m}+x^{15m}-x^{24m}+x^{27m}+x^{30m}-x^{33m}+x^{39m}+x^{45m}-x^{48m}+x^{51m}+x^{54m}+x^{57m}+x^{60m}.$
	 \end{itemize}
	 
	 It easy to see that $w_{2}(x)$ above satisfies 
	  $$ 3 \mid s \quad (\forall s \in P \cup N),$$ 
	  and that no polynomial equivalent to $w_1(x)$ has this property. 
	 Thus $ w'(x):=w_2(x^{\frac{1}{3}})$ is a solution in $CW(21,16)$:
	 $$ w'(x)=1-x^{1}-x^{2}-x^{4}+x^{5}-x^{8}+x^{9}+x^{10}-x^{11}+x^{13}+x^{15}-x^{16}+x^{17}+x^{18}+x^{19}+x^{20}.$$
	 Hence there is only one equivalence class in $ CW(21,16)$.  
	 This gives an equivalence class in $CW(21m,16)$ for each $m \geq 1$:

	 $ \widetilde{w'}(x)=1-x^{m}-x^{2m}-x^{4m}+x^{5m}-x^{8m}+x^{9m}+x^{10m}-x^{11m}+x^{13m}+x^{15m}-x^{16m}+x^{17m}+x^{18m}+x^{19m}+x^{20m}.$

{ \ }\\

\item  $ P=\{a,2a,4a,8a,b,2b,4b,8b,16b,32b\}; \  N=\{c,2c,d,2d,4d,8d\} $; \\  
	 $ ol(a)=4, \ ol(b)=6; \  ol(c)=2, \  ol(d)=4 $.\\

	 $P$ and $N$ contain orbits of lengths 2, 4 and  6.
	 By Section \ref{sec3}, $n$ must satisfy the conditions for the existence  
	 of orbits of length $i$, for $i \in \{2,4,6\}$: 
	 \begin{itemize}
	     \item \underline{$i=2$}: $n$ must be divisible by 3.
	     \item \underline{$i=4$}: $n$ must be divisible by 15; recall that if $5|n$ but 
		$15 \not| \; n$ then there is only one orbit of length 4.
	     \item \underline{$i=6$}: Here we need only one 
	     orbit of length 6 in $P \cup N$. Hence there are three possibilities: 
	     \begin{enumerate}
	\item $63 \mid n$, and then there are 9 different orbits of length 6.
	\item $21 \mid n$ but $63 \not| \; n$, and then there are 
	     two different orbits of length 6.    
	\item $9 \mid n$ but $63 \not| \; n $, and then there is only one 
	     orbit of length 6. 
	\end{enumerate}
	 
	 \end{itemize}
	 
	 Hence in the present case necessarily either $45 \mid n$ or $105 \mid n$.
	 The proofs of the following theorems are similar to those of Theorems  
	 \ref{th1} and \ref{th2}, and will be omitted.

\begin{thm} \label{th3}
{ \ }\\
(i)    For each odd $m \geq1 $, if $w(x) \in CW(45m,16)$ has (for the 
       multiplier $t=2$) 
		$$olp(P)=4^16^1, \quad olp(N)=2^14^1$$ 
       then:
       \begin{itemize}
	\item If $7 \not| \; m$ then there exists 
		$$w_0(x) \in CW(45,16) \quad {\rm s. \ t.} \quad w(x)=w_0(x^m).$$ 
	\item If $7 | m$ then there exists  
	$$w'_0(x) \in CW(315,16) \ {\rm s. \ t.} \ w(x)=w'_0(x^{\frac{m}{7}}).$$ 
	\end{itemize}
      
(ii)  If $w(x)$ and $\widetilde{w}(x)$ are equivalent in $R_{45m}$ then: 
	\begin{itemize}
	  \item 
		If  $7 \not| \; m$ then $w_0(x)$ and $\widetilde{w}_0(x)$ are 
		equivalent in $R_{45}$. 
	  \item 
		If $7 | m$ then $w'_0(x)$ and $\widetilde{w}'_0(x)$ are 
		equivalent in $R_{315}$. 
		
	\end{itemize}

\end{thm}
		   \begin{flushright}
				$\diamondsuit$
		   \end{flushright}

\begin{thm} \label{th4}
{ \ }\\
(i)    For each odd $m \geq1 $, if $w(x) \in CW(105m,16)$ has (for the 
       multiplier $t=2$) 
		$$olp(P)=4^16^1, \quad olp(N)=2^14^1$$ 
       then:
       \begin{itemize}
	\item If $3 \not| \; m$ then there exists 
		$$w_0(x) \in CW(105,16) \quad {\rm s. \ t.} \quad w(x)=w_0(x^m).$$ 
	\item If $3 | m$ then there exists  
	$$w'_0(x) \in CW(315,16) \ {\rm s. \ t.} \ w(x)=w'_0(x^{\frac{m}{3}}).$$ 
	\end{itemize}
      
(ii)  If $w(x)$ and $\widetilde{w}(x)$ are equivalent in $R_{105m}$ then: 
	\begin{itemize}
	  \item 
		If  $3 \not| \; m$ then $w_0(x)$ and $\widetilde{w}_0(x)$ are 
		equivalent in $R_{105}$. 
	  \item 
		If $3 | m$ then $w'_0(x)$ and $\widetilde{w}'_0(x)$ are 
		equivalent in $R_{315}$. 
		
	\end{itemize}

\end{thm}
    
		   \begin{flushright}
				$\diamondsuit$
		   \end{flushright}

	 By Theorems \ref{theor1}, \ref{th3} and \ref{th4} 
	 we have to search for circulant weighing only in $CW(315,16)$. 
	 This search will also give all solutions in $CW(45,16)$ and in 
	 $CW(105, 16)$. The data that we have consist of 
		$$n=315, \quad k=16, \quad t=2, \quad olp(P)=4^16^1, \quad olp(N)=2^14^1.$$ 
	 A Pascal program,
	 completely analogous to the ones
	 described in the previous cases, was written and run . No solutions 
	 were found. Hence there does not exist $w(x) \in CW(n,16)$ with the above data.

\end{enumerate}
	    
 %               }

\section{Summary} \label{sec7}

The following results were obtained in this paper. 
\begin{itemize}
    \item $CW(21,16) \neq \emptyset$. Only one equivalence class exists 
	  here (it was known before this work): 

	  $ w_{0}(x)=1-x^{1}-x^{2}-x^{4}+x^{5}-x^{8}+x^{9}+x^{10}-x^{11}+x^{13}+x^{15}-x^{16}+x^{17}+x^{18}+x^{19}+x^{20}.$

    \item $CW(31,16) \neq \emptyset$. There are two distinct equivalence 
	  classes (both were known before):
	  
	  $w_{1}(x)=-1-x-x^{2}+x^{3}-x^{4}+x^{6}+x^{7}-x^{8}+x^{12}+x^{14}-x^{16}+x^{17}+x^{19}+x^{24}+x^{25}+x^{28}, $
	  
	  $w_{2}(x)=-1-x-x^{2}-x^{4}+x^{5}-x^{8}+x^{9}+x^{10}+x^{15}-x^{16}+x^{18}+x^{20}+x^{23}+x^{27}+x^{29}+x^{30}. $

    \item $CW(63,16) \neq \emptyset$. Two distinct equivalence classes exist 
	  in this case: 
	  
	$w'_{1}(x)=1-x^{1}-x^{2}-x^{4}-x^{8}+x^{9}+x^{13}-x^{16}+x^{18}+x^{19}+x^{26}-x^{32}+x^{36}+x^{38}+x^{41}+x^{52}, $ 
\begin{flushright}
		[a new class which wasn't known before this work] 
\end{flushright}
	  
	$w'_{2}(x)=1-x^{3}-x^{6}-x^{12}+x^{15}-x^{24}+x^{27}+x^{30}-x^{33}+x^{39}+x^{45}-x^{48}+x^{51}+x^{54}+x^{57}+x^{60}.$
\begin{flushright}
		[an old class: $w'_{2}(x)=w_{0}(x^3)$ with $w_0(x) \in CW(21,16)$ above]
\end{flushright}   

	\item $CW(n,16) \neq \emptyset$, for odd n, iff either $21 | n$ or  
	      $31|n$. In each case, representatives for all possible 
	      equivalence classes are obtained by replacing $x$ by $x^m$ 
	      in one of the above polynomials, for a suitable integral value 
	      of $m$ ($ \frac{n}{21}, \quad \frac{n}{31}, \ {\rm or} \  \frac{n}{63}$). The number of equivalence classes is at most 4. 
	      More specifically: 
	      
	      If $31|n$ and $63|n$ then there are 4 classes.

	      If $31|n$ and $21|n$ but $63 \not| \; n$ then there are 3 classes. 

	      If $31|n$ but $21 \not| \; n$ then there are 2 classes.

	      If $31 \not| \; n$ but $63|n$ then there are 2 classes.

	      If $31 \not| \; n$ and $63 \not| \; n$ but $21 | n$ then there is one class. 

	      Otherwise ($31 \not| \; n$ and $21 \not| \; n$) - there are no classes.

\end{itemize}

%\chapter{}

}


\begin{thebibliography}{1}
  \bibitem{ams97}
    R.M. Adin, M. Muzychuk and Y. Strassler, Circulant weighing matrices of 
    prime order via symmetric polynomials, in preparation.
  
  \bibitem{adjp95}
    K.T. Arasu, J.F. Dillon, D. Jungnickel and A. Pott, The solution of the
    Waterloo problem, J. of Combinatorial Theory (Series A) 71, 316-331, 1995.
  
  \bibitem{as96}
    K.T. Arasu and J. Seberry, Circulant weighing designs, J. of 
    Combinatorial Designs 4, 439-447, 1996.
  
  \bibitem{as98}
    K.T. Arasu and J. Seberry, On circulant weighing matrices, preprint, 1998.  

%  \bibitem{as97}
%0    K.T. Arasu and Y. Strassler, A new class of circulant weighing matrices 
%    of order $24t$ with weight 9, in preparation.

  \bibitem{kal75}
    K.S. Banerjee, Weighing Designs for Chemistry, Medicine, Economics, Operations 
    Research, Statistics, M. Dekker, New York, 1975.
  

  
  \bibitem{cra91}
    R. Craigen, Constructions for Orthogonal Matrices, Ph.D. Thesis, 
    The  University of Waterloo, Ontario, Canada, 1991.

%  \bibitem{cra92}
%    R. Craigen, The structure of weighing matrices having large weights, 
%    Preprint.
  

  \bibitem{ead80}
    P. Eades, Circulant $(v,k, \mu)$ designs, Lecture Notes in Mathematics 
    No. 829, 83-93, Springer-Verlag, Berlin-Heidelberg, 1980.
  
  \bibitem{eh76}
    P. Eades, R.M. Hain, On circulant weighing matrices, Ars Combinatoria 2, 
    265-284, 1976.
  

  \bibitem{ggs76} 
    A.V. Geramita, J.M. Geramita, J. Seberry-Wallis, Orthogonal designs, 
    J. Lin. Mult. Algebra 3, 281-306, 1975.

  \bibitem{gs79}
    A.V. Geramita, J. Seberry. Orthogonal Designs (Quadratic Forms and
    Hadamard Matrices), Marcel Dekker, New York, 1979.

  \bibitem{hai77}
    R. Hain, Circulant Weighing Matrices, M.A. Thesis, Australian 
    National University, 1977.

  \bibitem{hal79}
    E.L. Hall, Computer Image Processing and Recognition, Academic Press, 
    New York, 1979.

  \bibitem{hal86}
    M. Hall, Combinatorial Theory, 2nd Ed., Wiley-Interscience, New York, 1986.

 \bibitem{hs79}
    M. Harmit and N.J.A. Sloan, Hadamard Transform Optics, Academic Press, 
    New York, 1979.

  \bibitem{jun92}
    D. Jungnickel, Difference sets, in: Jeffrey H. Dinitz and Douglas R. Stinton (editors),  
    Contemporary Design Theory: A Collection 
    of Surveys, John Wiley $\&$ Sons, 1992.

  \bibitem{koe94}
    G.J. Koehler, A proof of the Vose-Liepins conjecture, Annals of Mathematics 
    and Artificial Intelligence 10, 409-422, 1994.

%  \bibitem{ks96}
%    C. Koukouvinos and J. Seberry, New weighing matrices constructed using 
%    two sequences with zero autocorrelation function, preprint, 1996.

%  \bibitem{lss89}
%    C.W.H. Lam, T. Swiercz and S. Swiercz, The nonexistence of finite 
%    projective planes of order 10. Canad. J. Math. 41, 
%    1117-1123, 1989.
  
%  \bibitem{lam91}
%    C.W.H. Lam, The search for a finite projective plane of order 10, American 
%    Mathematical Monthly 98, 305-318, 1991.
  
  \bibitem{lan83}
    E.S. Lander, Symmetric Designs: An Algebraic Approach, Cambridge University 
    Press, 1983.
  

  \bibitem{mcf80}
    R.L. McFarland, On Multipliers of Abelian Difference Sets, 
    Ph.D. Thesis, Ohio State University, 1980.

%  \bibitem{ms97}
%    J.B. Muskat and Y. Strassler, Classification of integer weighing circulants 
%    of low orders up to 5 and any possible weight, in preparation.

  
  \bibitem{muz97}
    M. Muzychuk, Difference Sets with $n=2p^m$, J. Algebraic Combinatorics 7, 
    77-89, 1998. 


%  \bibitem{seb96}
%    J. Seberry, New asymptotic existence results for weighing matrices, 
%    preprint, 1996.
  
  \bibitem{sw75}
    J. Seberry-Wallis and A.L. Whiteman, Some results on weighing matrices, 
    Bull. Austral. Math. Soc. 12, 433-447, 1975.

  \bibitem{sw78}
    J. Seberry and K. Wehrhahn, A class of codes generated by circulant 
    weighing matrices, in:  D.H. Holton and J. Seberry (editors), 
    Proceedings of the International Conference on  
    Combinatorial Mathematics, Canberra, August 
    16-27, 1977, 
    282-289, Springer, Berlin, 1978.

  \bibitem{sy92}
    J. Seberry and M. Yamada, Hadamard matrices, sequences and block designs, 
    in: Jeffrey H. Dinitz  and Douglas R. Stinton (editors), Contemporary Design Theory: 
    A Collection of Surveys, John Wiley $\&$ Sons, 431-554, 1992.

  \bibitem{sm76}
    R.G. Stanton and R.C. Mullin. On existence of a class of circulant balanced 
    weighing matrices, SIAM J. Appl. Math. 30, 98-102, 1976.

  \bibitem{str83}
    Y. Strassler, In Search for Circulant Weighing Matrices, M.Sc. Thesis,
    Bar-Ilan University, 1983.

  \bibitem{str93}
     Y. Strassler, Circulant weighing matrices of prime order and weight 9 
     having a multiplier, manuscript.

  \bibitem{str95}
    Y. Strassler, New circulant weighing matrices of prime order in $CW(31,16)$, 
    $CW(71,25)$, $CW(127,64)$, Proceedings of the Bose Conference, J. qqStatistical 
    Planning and Inference 73, 317-330, 1998.

%  \bibitem{str97}
%    Y. Strassler, Quasisymmetric weighing circulants may exist only if their 
%    order and weight are even, in preparation.

  \bibitem{str}
    Y. Strassler, The Classification of Circulant Weighing Matrices of Weight 
    9, Ph.D. Thesis, Bar-Ilan University, 1998.


%  \bibitem{wil44}
%    J. Williamson, Hadamard's determinant theorem and the sum of four squares, 
%    Duke Math. J. 11, 65-81, 1944.


\end{thebibliography}
\end{document}